# GOODNESS-OF-FIT TESTS VIA PHI-DIVERGENCES

By Leah Jager[1] and Jon A. Wellner[2]

*Grinnell College and University of Washington*

A unified family of goodness-of-fit tests based on $\phi$-divergences is introduced and studied. The new family of test statistics $S_n(s)$ includes both the supremum version of the Anderson–Darling statistic and the test statistic of Berk and Jones [*Z. Wahrsch. Verw. Gebiete* **47** (1979) 47–59] as special cases ($s = 2$ and $s = 1$, resp.). We also introduce integral versions of the new statistics.

We show that the asymptotic null distribution theory of Berk and Jones [*Z. Wahrsch. Verw. Gebiete* **47** (1979) 47–59] and Wellner and Koltchinskii [*High Dimensional Probability III* (2003) 321–332. Birkhäuser, Basel] for the Berk–Jones statistic applies to the whole family of statistics $S_n(s)$ with $s \in [-1, 2]$. On the side of power behavior, we study the test statistics under fixed alternatives and give extensions of the "Poisson boundary" phenomena noted by Berk and Jones for their statistic. We also extend the results of Donoho and Jin [*Ann. Statist.* **32** (2004) 962–994] by showing that all our new tests for $s \in [-1, 2]$ have the same "optimal detection boundary" for normal shift mixture alternatives as Tukey's "higher-criticism" statistic and the Berk–Jones statistic.

**1. Introduction.** In this paper we introduce and study a new family of goodness-of-fit tests which includes both the supremum version of the Anderson–Darling statistic (or, equivalently, Tukey's "higher criticism" statistics as discussed by Donoho and Jin [15]) and the test statistic of Berk and Jones [5] as special cases. The new family is based on phi-divergences somewhat analogously to the phi-divergence tests for multinomial families introduced by Cressie and Read [10], and is indexed by a real parameter

Received March 2006; revised December 2006.

[1]Supported in part by NSF Grants DMS-05-03822 and DMS-03-54131 (VIGRE Grant).

[2]Supported in part by NSF Grants DMS-02-03320 and DMS-05-03822 and by NI-AID Grant 2R01 AI291968-04.

$s \in \mathbb{R}$: $s = 2$ gives the Anderson–Darling test statistic, $s = 1$ gives the Berk–Jones test statistic, $s = 1/2$ gives a new (Hellinger–distance type) statistic, $s = 0$ corresponds to the "reversed Berk–Jones" statistic studied by Jager and Wellner [24] and $s = -1$ gives a "Studentized" (or empirically weighted) version of the Anderson–Darling statistic. We introduce the corresponding integral versions of the new statistics (but will study them in detail elsewhere). Having a family of statistics available gives the possibility of better understanding of individual members of the family, as well as the ability to select particular members of the family that have different desirable properties.

In Section 2 we introduce the new test statistics. In Section 3 we briefly discuss the null distribution theory of the entire family of statistics, and note that the exact distributions can be handled exactly for sample sizes up to $n = 3000$ via Noé's recursion formulas (and possibly up to $n = 10^4$ via the recursion of Khmaladze and Shinjikashvili [31]) along the lines explored for the Berk–Jones statistic by Owen [36]. We also generalize the asymptotic distribution theory of Jaeschke [22] and Eicker [17] for the supremum version of the Anderson–Darling statistic, and of Berk and Jones [5] and Wellner and Koltchinskii [43] for the Berk–Jones statistic, by showing that the existing null distribution theory for $s = 1$ and $s = 2$ applies to (an appropriate version of) the whole family of statistics. We generalize the results of Owen [36] by showing that our family of test statistics provides a corresponding family of confidence bands.

In Section 4 we study the behavior of the new family of test statistics under fixed alternatives. We show that for $0 < s < 1$ and fixed alternatives the test statistics always converge almost surely to their corresponding natural parameters. For $1 < s < \infty$, we provide necessary and sufficient conditions on the alternative d.f. $F$ for convergence to the corresponding natural parameter to hold, and show that the "Poisson boundary" phenomena noted by Berk and Jones for their statistic continues to hold for $1 \le s < \infty$ and for $s < 0$ by identifying the Poisson boundary distributions explicitly. We also briefly discuss further large deviation results and connections between the work of Berk and Jones [5] and Groeneboom and Shorack [19].

In Section 5 we extend the results of Donoho and Jin [15] by showing that all our new tests for $s \in [-1, 2]$ have the same "optimal detection boundary" for normal shift mixture alternatives as Tukey's "higher-criticism" statistic and the Berk–Jones statistic.

Our new family of test statistics not only provides a unifying framework for the study of a number of existing test statistics as special cases, but also gives the possibility of "designing" a new test with several different desirable properties. For example, the new statistic $S_n(1/2)$ satisfies both: (a) it consistently estimates its "natural parameter" for every alternative and (b) it has the same optimal detection boundary for the two point normal mixture alternatives of Donoho and Jin [15] as the existing test statistics $S_n(2)$ and $S_n(1)$ considered by these authors.

**2. The test statistics.** Consider the classical goodness-of-fit problem: suppose that $X_1, \ldots, X_n$ are i.i.d. $F$, and let $\mathbb{F}_n(x) = n^{-1} \sum_{i=1}^n 1\{X_i \leq x\}$ be the empirical distribution function of the sample. We want to test

$$H : F = F_0 \quad \text{versus} \quad K : F \neq F_0,$$

where $F_0$ is continuous. By the probability integral transformation, we can, without loss of generality, suppose that $F_0$ is the uniform distribution on $[0,1]$, $F_0(x) = (x \wedge 1) \vee 0$, and that all the distribution functions $F$ in the alternative $K$ are defined on $[0,1]$. The basic idea behind our new family of tests is simple. For fixed $x \in (0,1)$, the interval is divided into two sub-intervals $[0,x]$ and $(x,1]$, and we can test the (pointwise) null hypothesis $H_x : F(x) = x$ versus the (pointwise) alternative $K_x : F(x) \neq x$ using any of the general phi-divergence test statistics $K_\phi(\mathbb{F}_n(x), x)$ proposed by Csiszár [11] (see also Csiszár [12] and Ali and Silvey [2]) and studied further in a multinomial context by Cressie and Read [10], where $\phi$ is a convex function mapping $[0, \infty)$ to the extended reals $\mathbb{R} \cup \{\infty\}$ (cf. Liese and Vajda [32], pages 10 and 212, and Vajda [39]). Then our proposed test statistics are of the form

$$S_n(\phi) \equiv \sup_x K_\phi(\mathbb{F}_n(x), x)$$

or

$$T_n(\phi) \equiv \int_0^1 K_\phi(\mathbb{F}_n(x), x) \, dx,$$

where the supremum and/or integral over $x$ may require some restriction depending on the choice of $\phi$.

In our particular case, we define $\phi = \phi_s$ for $s \in \mathbb{R}$ by

$$\phi_s(x) \equiv \begin{cases} [1 - s + sx - x^s]/[s(1-s)], & s \neq 0, 1, \\ x(\log x - 1) + 1 \equiv h(x), & s = 1, \\ \log(1/x) + x - 1 \equiv \tilde{h}(x), & s = 0 \end{cases}$$

(cf. Liese and Vajda [32], page 34), so that

$$: K_s(u, v) = v\phi_s(u/v) + (1-v)\phi_s((1-u)/(1-v))$$
$$= \frac{1}{s(1-s)}\{1 - u^s v^{1-s} - (1-u)^s(1-v)^{1-s}\}, \qquad s \neq 0, 1.$$

Note that this definition makes $\phi_s$ continuous in $s$ for all $x$ in $(0,1)$, and hence, $K_s$ is continuous in $s$ for all $(u, v) \in (0,1)^2$. Also note that $K_{\lambda+1}(p, q) = I_2^\lambda(\underline{p} : \underline{q})$, where $\underline{p} = (p, 1-p)$, $\underline{q} = (q, 1-q)$, and $I_2^\lambda(\underline{p} : \underline{q})$ is as defined in (5.1),

Cressie and Read [10], page 456. Then our proposed test statistics $S_n(s)$ and $T_n(s)$ for $s \in \mathbb{R}$ are defined by

$$(1) \qquad S_n(s) \equiv \begin{cases} \sup_{0 < x < 1} K_s(\mathbb{F}_n(x), x), & \text{if } s \geq 1, \\ \sup_{X_{(1)} \leq x < X_{(n)}} K_s(\mathbb{F}_n(x), x), & \text{if } s < 1 \end{cases}$$

and

$$(2) \qquad T_n(s) \equiv \begin{cases} \displaystyle\int_0^1 K_s(\mathbb{F}_n(x), x)\,dx, & \text{if } s > 0, \\ \displaystyle\int_{X_{(1)}}^{X_{(n)}} K_s(\mathbb{F}_n(x), x)\,dx, & \text{if } s \leq 0. \end{cases}$$

The reasons for changing the definitions of the statistics by restricting the supremum or integral for different values of $s$ will be explained in Section 3; basically, the restrictions must be imposed for some appropriate value of $s$ in order to maintain the same null distribution theory for all values of $s$ in $[-1, 2]$.

The most notable special cases of these statistics are $s \in \{-1, 0, 1/2, 1, 2\}$: it is easily checked that

$$K_2(u, v) = \frac{1}{2}\frac{(u - v)^2}{v(1 - v)},$$

$$K_1(u, v) = u \log\left(\frac{u}{v}\right) + (1 - u) \log\left(\frac{1 - u}{1 - v}\right),$$

$$K_{1/2}(u, v) = 4\{1 - \sqrt{uv} - \sqrt{(1 - u)(1 - v)}\}$$
$$= 2\{(\sqrt{u} - \sqrt{v})^2 + (\sqrt{1 - u} - \sqrt{1 - v})^2\},$$

$$K_0(u, v) = K_1(v, u) = v \log\left(\frac{v}{u}\right) + (1 - v) \log\left(\frac{1 - v}{1 - u}\right),$$

$$K_{-1}(u, v) = K_2(v, u) = \frac{1}{2}\frac{(u - v)^2}{u(1 - u)}.$$

It follows that:

(a) $S_n(2)$ is (1/2 times) the square of the supremum form of the Anderson–Darling statistic (or, in its one-sided form, Tukey's "higher criticism statistic"; see Donoho and Jin [15] and Section 5).

(b) $S_n(1)$ is the statistic studied by Berk and Jones [5].

(c) $S_n(1/2)$ is (4 times) the supremum of the pointwise Hellinger divergences between Bernoulli($\mathbb{F}_n(x)$) and Bernoulli($F_0(x)$); as far as we know, this is a new goodness-of-fit statistic [as are all the statistics $S_n(s)$ for $s \notin \{-1, 0, 1, 2\}$].

(d) $S_n(0)$ is the "reversed Berk–Jones" statistic introduced by Jager and Wellner [24].

(e) $S_n(-1)$ is (1/2 times) a "Studentized" version of the supremum form of the Anderson–Darling statistic; see, for example, Eicker [17], page 116.

(f) $T_n(1)$ is the integral form of the Berk–Jones statistic introduced by Einmahl and McKeague [18].

(g) $T_n(2)$ is the classical (integral form of) the Anderson–Darling statistic introduced by Anderson and Darling [3].

REMARK 2.1. Note that $K_{1/2-r}(u,v) = K_{1/2+r}(v,u)$ for $r \in \mathbb{R}$ and $u$, $v \in (0,1)$, so the families of statistics $S_n(s)$ and $T_n(s)$ have a natural symmetry about $s = 1/2$. We will continue to use the "$s$-parametrization" of these families for reasons of notational simplicity.

## 3. Distributions under the null hypothesis.

3.1. *Finite sample critical points via Noé's recursion.* Owen [36] showed how to use the recursions of Noé [35] to obtain finite sample critical points of the Berk–Jones statistic $R_n = S_n(1)$ for values of $n$ up to 1000. (See Shorack and Wellner [38], pages 362–366 for an exposition of Noé's methods.) Jager and Wellner [24] pointed out a minor error in the derivations of Owen [36] and extended his results to the reversed Berk–Jones statistic $S_n(0)$. Jager [23] gives exact finite sample computations for the whole family of statistics via Noé's recursions for values of $n$ up to 3000. (The C and R programs are available at the second author's website.) We will not give details of the finite-sample computations here but refer the interested reader to Jager and Wellner [24] and Jager [23]. See Jager and Wellner [24] and Jager [23] for plots of finite sample critical points and several finite sample approximations based on the asymptotic theory given here.

During the revision of this paper we learned of an alternative finite-sample recursion for calculation of the null distribution of $S_n(2)$ proposed by Khmaladze and Shinjikashvili [31] which apparently works for $n \leq 10^4$. Presumably this alternative recursion could be used for our entire family of statistics, but this has not yet been carried out.

3.2. *Asymptotic distribution theory for $S_n(s)$ under the null hypothesis.* Limit distribution theory for $S_n(2)$ and $S_n(-1)$ under the null hypothesis follows from the work of Jaeschke [22] and Eicker [17]; see Shorack and Wellner [38], Chapter 16, pages 597–615 for an exposition. These results are closely related to the classical results of Darling and Erdös [14]. Berk and Jones [5] stated the asymptotic distribution of their statistic $R_n = S_n(1)$. For details of the proof, see Wellner and Koltchinskii [43], with a minor correction as noted at the end of the proof here. Here we show that the limit distribution of $nS_n(s) - r_n$ is the same double-exponential extreme value distribution for all $-1 \leq s \leq 2$, where

$$r_n = \log_2 n + \tfrac{1}{2} \log_3 n - \tfrac{1}{2} \log(4\pi),$$

with $\log_2 n \equiv \log(\log n)$ and $\log_3 n \equiv \log(\log_2 n)$.

THEOREM 3.1 (Limit distribution under null hypothesis). *Suppose that the null hypothesis $H$ holds so that $F$ is the uniform distribution on $[0,1]$. Then for $-1 \leq s \leq 2$ it follows that*

$$nS_n(s) - r_n \xrightarrow{d} Y_4 \sim E_v^4,$$

*where $E_v^4(x) = \exp(-4\exp(-x)) = P(Y_4 \leq x)$.*

Define

$$b_n = \sqrt{2\log_2 n}, \qquad c_n = 2\log_2 n + \tfrac{1}{2}\log_3 n - \tfrac{1}{2}\log(4\pi),$$

$$d_n = n^{-1}(\log n)^5, \qquad Z_n \equiv \sup_{d_n \leq x \leq 1 - d_n} \frac{\sqrt{n}|\mathbb{F}_n(x) - x|}{\sqrt{x(1-x)}}.$$

As will be seen, the proof involves the following four facts: *Fact* 1. $Z_n/b_n \to_p 1$. *Fact* 2. $b_n Z_n - c_n \xrightarrow{d} Y_4 \sim E_v^4$. *Fact* 3. $(1/2)c_n^2/b_n^2 = r_n + o(1)$. *Fact* 4. $nS_n(s) = (1/2)Z_n^2 + o_p(1)$.

In the ranges $s > 2$ and $s < -1$ we do not know a theorem describing the behavior of the statistics $S_n(s)$ under the null hypothesis.

3.3. *Confidence bands.* Owen [36] showed how the Berk–Jones statistic $R_n = S_n(1)$ can be inverted to obtain confidence bands for an unknown distribution function $F$. Similarly, the family of statistics $S_n(s)$ yields a new family of confidence bands for $F$ as follows: given a continuous d.f. $F$ on $\mathbb{R}$, define

$$S_n(s, F) \equiv \begin{cases} \sup_{-\infty < x < \infty} K_s(\mathbb{F}_n(x), F(x)), & \text{if } s \geq 1, \\ \sup_{X_{(1)} \leq x < X_{(n)}} K_s(\mathbb{F}_n(x), F(x)), & \text{if } s < 1. \end{cases}$$

By the (inverse) probability integral transformation, $P_F(S_n(s, F) \leq t) = P_{F_0}(S_n(s) \leq t)$ for all $t$ where $S_n(s)$ is as defined in (1) and $F_0$ is the uniform distribution on $[0,1]$. Hence, with $q_n(s, \alpha)$ denoting the upper $1 - \alpha$ quantile of the distribution of $S_n(s)$ under $F_0$ (which is computable via Noé's recursion as discussed in Section 3.1 or can be approximated for large $n$ via Theorem 3.1), it follows that

$$P_F(S_n(s, F) \leq q_n(s, \alpha)) = P_{F_0}(S_n(s) \leq q_n(s, \alpha)) = 1 - \alpha$$

for each fixed $\alpha \in (0, 1)$ and $n$. Hence,

$$\{F : S_n(s, F) \leq q_n(s, \alpha)\} = \{F : L_n(x; s, \alpha) \leq F(x) \leq U_n(x; s, \alpha) \text{ for all } x \in \mathbb{R}\}$$

yields a family of $1 - \alpha$ confidence bands for $F$. Here $L_n(x; s, \alpha)$ and $U_n(x; s, \alpha)$ are random functions determined by $s$, $\alpha$, $n$ and the data in a straightforward way; see Owen [36], Jager and Wellner [24] and Jager [23] for details.

3.4. *Asymptotic distribution theory for $T_n(s)$ under the null hypothesis.* Limit distribution theory for $T_n(2)$ was established by Anderson and Darling [3]. Einmahl and McKeague [18] noted that this carries over to $T_n(1)$ (for a proof, see Wellner and Koltchinskii [43]) and extended $T_n(1)$ to other testing problems. Here we show that the limit distribution of $nT_n(s)$ is ($1/2$ times) the Anderson–Darling limit distribution for all $s \in [-1, 2]$, namely, the distribution of

$$(3) \qquad A^2 \equiv \int_0^1 \frac{[\mathbb{U}(t)]^2}{t(1-t)} \, dt \overset{d}{=} \sum_{j=1}^{\infty} \frac{Z_j^2}{j(j+1)},$$

where $\mathbb{U}$ is a standard Brownian bridge process on $[0, 1]$ and $Z_1, Z_2, \ldots$ are i.i.d. $N(0, 1)$; see, for example, Shorack and Wellner [38], pages 224–227.

THEOREM 3.2 (Limit distribution of $T_n(s)$ under the null hypothesis). *Suppose that the null hypothesis $H$ holds so that $F$ is the uniform distribution on $[0, 1]$. Then for $-\infty < s \leq 2$ it follows that $nT_n(s) \overset{d}{\to} A^2/2$, where $A^2$ is the Anderson–Darling limit defined in* (3).

We will not study the statistics $T_n(s)$ further in this paper, but intend to continue their study elsewhere.

**4. Limit theory under alternatives.** The power behavior of individual members of our new family of statistics has previously been studied separately and somewhat in isolation: see, for example, Berk and Jones [5] (for the Berk–Jones statistic), Durbin, Knott and Taylor [16] and D'Agostino and Stephens [13] for $T_n(2)$ compared to other integral goodness-of-fit statistics and Nikitin [34] for treatment of Bahadur efficiencies for many goodness-of-fit statistics. Interest in these test statistics has received new impetus via the use of appropriate one-sided versions of the test statistic $S_n(2)$ in the context of multiple testing problems; see, for example, Donoho and Jin [15], Jin [28] and Meinshausen and Rice [33]. See Cayón, Jin and Treaster [8] for an interesting application to detection of non-Gaussianity in the cosmic microwave background data gathered by the Wilkinson Microwave Anisotropy Probe (WMAP) satellite, and see Cai, Jin and Low [7] for further work on estimation aspects of the problem in connection with the developments in Meinshausen and Rice [33]. The work of Owen [36] on confidence bands derived from the Berk–Jones statistic $S_n(1)$ was apparently motivated in large part by the Bahadur efficiency results of Berk and Jones [4].

It is clear from the results of Donoho and Jin [15] and earlier efforts by Révész [37] to combine the strengths of the Kolmogorov and Jaeschke–Eicker statistics that tests based on any of the statistics $S_n(s)$ will do best against "alternatives in the tails." As suggested by one of the referees of our paper, this may well be the "Achilles heel" of such test statistics, since the results of Révész [37] suggest that our statistics will have no asymptotic power for a large class of "contiguous alternatives" (with departures from the null hypothesis "in the middle" of the distribution). On the other hand, having a family of statistics such as $\{S_n(s): s \in \mathbb{R}\}$ available gives the possibility of choosing (or "designing") a test with several desirable properties. We will return to this briefly in Section 6.

Here we study convergence of the family of statistics to their "natural parameters" under fixed alternatives, comment briefly on the Bahadur efficiency results of Berk and Jones [5] in light of the results of Groeneboom and Shorack [19], and show that the optimal detection boundary results of Donoho and Jin [15] extend to the whole family of statistics $S_n(s)$ for $s \in [-1, 2]$. In spite of the negative results of Janssen [27] for goodness-of-fit statistics in general, much remains to be learned about the power behavior of the family $\{S_n(s)\}$.

4.1. *Almost sure convergence to natural parameter.* Let $F_0$ be the Uniform$(0, 1)$ distribution function as in Section 2. The Kolmogorov statistic $D_n \equiv \|\mathbb{F}_n - F_0\|_\infty$ has the property that for *any* distribution function $F$, if $X_1, \ldots, X_n$ are i.i.d. $F$, then

$$D_n \overset{a.s.}{\to} \|F - F_0\|_\infty \equiv d(F).$$

We call $d(F) = \|F - F_0\|_\infty$ the *natural parameter* for the Kolmogorov statistic $D_n$. As Berk and Jones [5] pointed out for their statistic $R_n = S_n(1)$, under alternatives $F$ the convergence

$$S_n(1) = R_n = \sup_{0 < x < 1} K_1(\mathbb{F}_n(x), x) \overset{a.s.}{\to} \sup_{0 < x < 1} K_1(F(x), x) \equiv r(F)$$

holds only under some condition on $F$ (the exact condition will be given below), and for a slightly more extreme $F$, namely, what we call the "Poisson boundary distribution function," the behavior changes to convergence in distribution to a functional of a Poisson process rather than convergence to a natural parameter. Thus, Berk and Jones [5] showed that if $F(x) = 1/(1 + \log(1/x))$, then

$$(4) \qquad S_n(1) = R_n \overset{d}{\to} \sup_{t > 0} \frac{\mathbb{N}(t)}{t} \overset{d}{=} \frac{1}{U},$$

where $\mathbb{N}$ is a standard Poisson process and $U \sim \text{Uniform}(0,1)$.

It turns out that in the range $0 < s < 1$ the statistics $S_n(s)$ behave analogously to the Kolmogorov statistic $D_n$ under fixed alternatives. Namely, we show that in this range the statistics converge almost surely to their "natural parameter" for all d.f.'s $F$.

PROPOSITION 4.1. *Suppose that* $X_1, \ldots, X_n$ *are i.i.d.* $F$ *and that* $0 < s < 1$. *Then* $S_n(s) \overset{a.s.}{\to} \sup_{0 < x < 1} K_s(F(x), x) \equiv S_\infty(s, F)$.

On the other hand, in the range $s > 1$ we have the following criterion for almost sure convergence of the statistics $S_n(s)$ to their natural parameters:

PROPOSITION 4.2. *Suppose that* $X_1, \ldots, X_n$ *are i.i.d.* $F$ *and that* $s > 1$. *Then* $S_n(s) \overset{a.s.}{\to} \sup_{0 < x < 1} K_s(F(x), x) \equiv S_\infty(s, F)$ *if and only if* $F$ *satisfies*

$$\int_0^1 \frac{1}{(F^{-1}(u)(1 - F^{-1}(u)))^{(s-1)/s}} \, du < \infty.$$

*By the (inverse) probability integral transformation, the convergence in the last display is equivalent to* $E_F[X(1-X)]^{-(1-1/s)} < \infty$.

As Berk and Jones [5] show, if for some $\gamma > 0$, the distribution function $F$ satisfies

$$F(x) \le \{\log(1/x)(\log_2(1/x))^{1+\gamma}\}^{-1}, \qquad x \le \gamma,$$

and

$$1 - F(x) \le \{\log(1/(1-x))(\log_2(1/(1-x)))^{1+\gamma}\}^{-1}, \qquad x \ge 1 - \gamma,$$

then $R_n \equiv S_n(1) \overset{a.s.}{\to} \sup_{0 < x < 1} K_1(F(x), x) \equiv S_\infty(1, F) \equiv r(F) < \infty$. It can be shown that this convergence holds if and only if $\int_0^1 [x(1-x)]^{-1} F(x) \times (1 - F(x)) \, dx < \infty$; see Jager [23] for details. We do not yet know sharp conditions for $S_n(s) \overset{a.s.}{\to} S_\infty(s, F)$ when $s \le 0$.

**4.2. Poisson boundaries for $s \ge 1$ and $s < 0$.** As noted in the previous subsection, the statistic $R_n = S_n(1)$ has a "Poisson boundary" d.f. $F_1$ for which $R_n = S_n(1) \overset{d}{\to} 1/U$ rather than $R_n = S_n(1) \overset{a.s.}{\to} r(F) \equiv S_\infty(1, F)$. Here we note that this behavior persists for the entire range $s \ge 1$ and for $s < 0$.

For each fixed $s \in [0,1)^c$, define the distribution function $F_s$ on $[0,1]$ by $F_s(0) = 0$ and, for $0 < x \le 1$, by

$$(5) \qquad F_s(x) = \begin{cases} \left(1 + \dfrac{x^{1-s} - 1}{s - 1}\right)^{-1/s}, & 1 < s < \infty, \\ (1 + \log(1/x))^{-1}, & s = 1, \\ (1 - s(x^{s-1} - 1))^{1/s}, & s < 0. \end{cases}$$

Note that $F_s(x) \to F_1(x)$ as $s \searrow 1$ for $0 \le x \le 1$.

The following proposition includes the result (4) of Berk and Jones [5] when $s = 1$, and it agrees with the case $b = 1/2$ of Theorem 2 of Jager and Wellner [25] when $s = 2$. As in (4), let $\mathbb{N}$ be a standard Poisson process, and let $U \sim \text{Uniform}(0, 1)$.

PROPOSITION 4.3 (Poisson boundaries for $s \ge 1$ and $s < 0$).

(i) *Fix $s \ge 1$ and suppose that $X_1, \ldots, X_n$ are i.i.d. $F_s$ given in (5). Then*

$$S_n(s) \xrightarrow{d} \frac{1}{s}\left(\sup_{t > 0} \frac{\mathbb{N}(t)}{t}\right)^s \stackrel{d}{=} \frac{1}{sU^s}.$$

(ii) *Fix $s < 0$ and suppose that $X_1, \ldots, X_n$ are i.i.d. $F_s$ given in (5). Then*

$$S_n(s) \xrightarrow{d} \frac{1}{(1-s)}\left(\sup_{t \ge S_1} \frac{t}{\mathbb{N}(t)}\right)^{-s}, \tag{6}$$

*where $S_1 = E_1$ is the first jump point of $\mathbb{N}$.*

REMARK 4.1. The distribution of $\sup_{t \ge S_1}(t/\mathbb{N}(t))$, which is also the limiting distribution of $\sup\{(t/\mathbb{G}_n(t)) : t \ge \xi_{(1)}\}$ where $\mathbb{G}_n$ is the empirical distribution function of $n$ i.i.d. Uniform$(0, 1)$ random variables $\xi_1, \ldots, \xi_n$, is given by

$$P\left(\sup_{t \ge S_1}(t/\mathbb{N}(t)) > x\right) = \exp(-x) + \sum_{k=1}^{\infty} \frac{(k-1)^{k-1}}{k!} x^k \exp(-kx), \qquad x > 1;$$

see Wellner [40], pages 1008–1009 and Shorack and Wellner [38], page 412. This yields an explicit formula for the distribution of the random variable on the right-hand side of (6).

REMARK 4.2. Although the family of distributions $F_s$ satisfies $F_s(x) \to \exp(-(1/x - 1)) \equiv \widetilde{F}_0(x)$ as $s \nearrow 0$, it appears that the natural limit in distribution under $\widetilde{F}_0$ is $S_n(0) \xrightarrow{d} 1 = \sup_{0 < x < 1} K_0(\widetilde{F}_0(x), x)$ in this case, so apparently convergence to the natural parameter continues to hold under $\widetilde{F}_0$. We do not know if there is a (more extreme) d.f. $F_0^{\dagger}$ for which $S_n(0) \xrightarrow{d} g(\mathbb{N})$ for a nondegenerate functional $g$ of a standard Poisson process $\mathbb{N}$.

REMARK 4.3. If $F \in K$ has Poisson boundary behavior at both 0 and 1, then natural generalizations of Proposition 4.3 involving two independent Poisson processes can easily be proved. For example, if $F$ is the standard arcsin law with density $\pi^{-1} u^{-1/2}(1 - u)^{-1/2} 1_{(0,1)}(u)$, then $S_n(s) \xrightarrow{d} (2/\pi^2)\max\{(\sup_{t>0}(\mathbb{N}(t)/t))^2, (\sup_{t>0}(\widetilde{\mathbb{N}}(t)/t))^2\}$, where $\mathbb{N}$, $\widetilde{\mathbb{N}}$ are independent standard Poisson processes.

4.3. *Bahadur efficiency comparisons.* Berk and Jones [5] studied the Bahadur efficiency of their statistic $S_n(1) = R_n$ relative to weighted Kolmogorov statistics based on the work of Abrahamson [1]. As pointed out by Groeneboom and Shorack [19], however, the Bahadur efficacies of the weighted Kolmogorov statistics are 0 for weights heavier than the (quite light) logarithmic weight function $\psi(x) = -\log(x(1-x))$ because the null distribution large-deviation result is degenerate for heavier weights. A second difficulty for Bahadur efficiency comparisons is that both the weighted Kolmogorov statistics and the Berk–Jones statistic fail to converge almost surely to their natural parameters for sufficiently extreme alternative d.f.'s $F$, and as noted by Berk and Jones [5] for the Berk–Jones statistic and by Jager and Wellner [25] for the weighted Kolmogorov statistics, there is a certain "Poisson boundary" d.f. $F$ for which the statistics converge in distribution to a functional of a Poisson process. Thus, comparisons of goodness-of-fit statistics of the supremum type via Bahadur efficiency are rendered difficult by breakdowns in both the large-deviation theory under the null hypothesis and by failure of the statistics to converge almost surely under fixed alternatives. Nevertheless, it would be interesting to be able to make comparisons where possible.

To this end, we consider variants of our statistics in the range $0 < s < 1$ with the supremum unrestricted as follows:

$$S_n^{ur}(s) = \sup_{0 < x < 1} K_s(\mathbb{F}_n(x), x),$$

$$S_n^{ur,+}(s) = \sup_{0 < x < 1} K_s^+(\mathbb{F}_n(x), x), \qquad S_n^{ur,-}(s) = \sup_{0 < x < 1} K_s^-(\mathbb{F}_n(x), x),$$

where

$$K_s^+(u, v) = \begin{cases} K_s(u, v), & \text{if } 0 < v < u < 1, \\ 0, & \text{if } 0 \le u \le v \le 1, \\ \infty, & \text{otherwise,} \end{cases}$$

$$K_s^-(u, v) = \begin{cases} K_s(u, v), & \text{if } 0 < u < v < 1, \\ 0, & \text{if } 0 \le v \le u \le 1, \\ \infty, & \text{otherwise.} \end{cases}$$

Also, set $K(u, v) = K_1(u, v) = u \log(u/v) + (1 - u \log((1-u)/(1-v))$ for $(u, v) \in (0, 1)^2$ and $K^+(u, v) = K_1^+(u, v)$. Although we do not yet have large deviation results for the statistics $S_n(s)$ or $S_n^+(s) \equiv \sup_{X_{(1)} \le x < X_{(n)}} K_s^+(\mathbb{F}_n(x), x)$, we can establish the following large deviation results for $S_n^{ur,+}(s)$ and $S_n^{ur,-}(s)$.

THEOREM 4.4. *Suppose that $X_1, \ldots, X_n$ are i.i.d. with continuous d.f. $F_0$, the uniform distribution on $(0, 1)$. Fix $s \in (0, 1)$. Then*

$$n^{-1} \log P_0(S_n^{ur,+}(s) \ge a) \to -\inf_{0 < x < 1} K^+(\tau_s^+(x, a), x)$$

(7)
$$= -\frac{\log[1 - s(1 - s)a]}{1 - s} \equiv -g_s^+(a)$$

for each $0 \le a < 1/[s(1 - s)]$, where $\tau_s^+(x, a) = \inf\{t : K_s^+(t, x) \ge a\}$. Furthermore,

$$n^{-1} \log P_0(S_n^{ur, -}(s) \ge a) \to - \inf_{0 < x < 1} K^-(\tau_s^-(x, a), x)$$

$$= -\frac{\log[1 - s(1 - s)a]}{1 - s} \equiv -g_s^-(a)$$

for each $a \ge 0$, where $\tau_s^-(x, a) = \sup\{t : K_s^-(t, x) \ge a\}$.

Combining Theorem 4.4 with Proposition 4.1, we have the following corollary for the Bahadur efficacies of the statistics $S_n^{ur, +}(s)$ and $S_n^{ur, -}(s)$ with $0 < s < 1$.

COROLLARY 4.5. *Let $F$ be a continuous distribution function on $[0, 1]$. Then the Bahadur efficacy of $S^{ur, +}(s)$ at the alternative $F$ is*

$$\epsilon_s^\pm(F) = g_s^\pm(S_\infty^\pm(s, F)) = g_s^+(S_\infty^\pm(s, F)),$$

*where $g_s^+$ is defined in (7) and $S_\infty^\pm(s, F) \equiv \sup_{0 < x < 1} K_s^\pm(F(x), x)$.*

REMARK. Note that $\lim_{s \nearrow 1} g_s^+(a) = a$, in agreement with Theorem 2.2, page 50, of Berk and Jones [5].

REMARK. Since $g_s^+(a) = g_s^-(a) \sim sa$ as $s \searrow 0$, the Bahadur efficacies of the statistics $S_n^{ur, \pm}(s)$ tend to be smaller than the efficacies of the Berk–Jones statistic $S_\infty^+(1, F) = r(F)$ (when the latter exists), and especially so for small $s$. This, together with extensive numerical computations of Jager [23], strengthens the case in favor of the statistics $S_n^+(s) = \sup_{X_{(1)} \le x < X_{(n)}} K_s^+ \times (\mathbb{F}_n(x), x)$ with restricted supremum. Unfortunately we do not yet know the large deviation behavior of these statistics with restricted supremum.

## 5. Attainment of the Ingster–Donoho–Jin optimal detection boundary.
Jin [28] and Donoho and Jin [15] consider testing in a "sparse heterogenous mixture" problem defined as follows: Suppose that $Y_1, \ldots, Y_n$ are i.i.d. $G$ on $\mathbb{R}$ and consider testing

$$H_0 : G = \Phi, \qquad \text{the standard } N(0, 1) \text{ distribution function}$$

versus

$$H_1 : G = (1 - \epsilon)\Phi + \epsilon\Phi(\cdot - \mu) \qquad \text{for some } \epsilon \in (0, 1), \mu > 0.$$

In particular, they consider the $n$-dependent alternatives $H_1^{(n)}$ given by

$$(8) \quad H_1^{(n)} \colon G_n = (1 - \epsilon_n)\Phi + \epsilon_n \Phi(\cdot - \mu_n) \qquad \text{for } \epsilon_n = n^{-\beta}, \mu_n = \sqrt{2r \log n},$$

where $1/2 < \beta < 1$ and $0 < r < 1$. By transforming to $X_i \equiv 1 - \Phi(Y_i)$ i.i.d. $F = 1 - G(\Phi^{-1}(1 - \cdot))$ (with the $X_i$'s taking values in $[0, 1]$), the testing problem becomes test

$$H_0 \colon F = F_0, \qquad \text{the Uniform}(0, 1) \text{ distribution function}$$

versus

$$H_1 \colon F = F_0(u) + \epsilon\{(1 - u) - \Phi(\Phi^{-1}(1 - u) - \mu)\} > F_0(u).$$

[The corresponding $n$-dependent sequence is $F_n(u) = u + \epsilon_n\{(1 - u) - \Phi(\Phi^{-1}(1 - u) - \mu_n)\}$ with the same choice of $\epsilon_n$ and $\mu_n$ as in (8).] Donoho and Jin [15] consider several different test statistics, among which the principal contenders are Tukey's "higher criticism" statistic $HC_n^*$ defined by

$$HC_n^* \equiv \sup_{X_{(1)} \leq x < X_{([\alpha_0 n])}} \frac{\sqrt{n}(\mathbb{F}_n(x) - x)}{\sqrt{x(1 - x)}}$$

for some $\alpha_0 > 0$ (they seem to usually take $\alpha_0 = 1/2$), and a one-sided version of the Berk–Jones statistic $BJ_n^+ \equiv n \sup_{X_{(1)} \leq x < 1/2} K_1^+(\mathbb{F}_n(x), x)$, where $K_s^+(u, v) \equiv K_s(u, v) 1\{0 < v < u < 1\}$.

Jin [28] (see also Ingster [20, 21]) showed that the likelihood ratio test of $H_0$ versus $H_1^{(n)}$ has a "detection boundary" defined in terms of the parameters $\beta \in (1/2, 1)$ and $r \in (0, 1)$ involved in (8) which is described as follows: set

$$\rho^*(\beta) = \begin{cases} \beta - 1/2, & 1/2 < \beta \leq 3/4, \\ (1 - \sqrt{1 - \beta})^2, & 3/4 < \beta < 1. \end{cases}$$

Then for $r > \rho^*(\beta)$, the likelihood ratio test (which makes use of knowledge of $\beta$ and $r$) is size and power consistent against $H_1^{(n)}$ as $n \to \infty$. Donoho and Jin [15] show that the tests of $H_0$ versus $H_1^{(n)}$ based on $HC_n^*$ and $BJ_n^+$ are also size and power consistent as $n \to \infty$ and that both of these tests dominate several other tests based on multiple comparison procedures such as the sample range, sample maximum, FDR (False Discovery Rate) and Fisher's method; see, for example, Figure 1 of Donoho and Jin [15] and their Theorems 1.4 and 1.5.

We show here that the tests based on appropriate one-sided versions of the statistics $S_n(s)$, namely,

$$nS_n^+(s) \equiv n \sup_{X_{(1)} \leq x \leq 1/2} K_s^+(\mathbb{F}_n(x), x),$$

have the same detection boundary for testing $H_0$ versus $H_1^{(n)}$ as the statistics $HC_n^*$ and $BJ_n^+$. More formally, define a function $\rho_s(\beta)$ such that if $\mu_n = \sqrt{2r \log n}$ and if we use a sequence of levels $\alpha_n \to 0$ slowly enough [slowly enough so that with $q_n(s, \alpha)$ as defined in Section 3.2, $\alpha_n$ satisfies $nq_n(s, \alpha_n) = (1 + o(1)) \log \log n$], then for $r > \rho_s(\beta)$, the resulting sequence of tests has power tending to 1 as $n \to \infty$, while for $r < \rho_s(\beta)$, the sequence of tests has power tending to zero. In terms of the functions $\rho_s(\beta)$, our theorem is as follows.

THEOREM 5.1.   *For each $s \in [-1, 2]$, $\rho_s(\beta) = \rho^*(\beta)$ for $1/2 < \beta < 1$.*

While this may not be too suprising for $1 \le s \le 2$ in view of the Donoho–Jin results for $nS_n^+(1) = BJ_n^+$ and $nS_n^+(2) = (1/2)(HC_n^*)^2$, it seems new and interesting for $s \in [-1, 1)$. Figure 1 gives smoothed histograms of the values of the statistics $nS_n(s) - r_n$ under the null hypothesis $H_0$ (solid line) and under the alternative hypothesis $H_1^{(n)}$ (dotted line) for $n = 0.5 \times 10^6$, $r = 0.15$ and $\beta = 1/2$. This should be compared with Figure 2 on page 978 of Donoho and Jin [15] showing values of $HC_n^*$ and $HC_n^+$, corresponding to our $s = 2$; their $HC_n^+ \equiv \sup_{1/n \le x \le 1/2} \sqrt{n} (\mathbb{F}_n(x) - x)^+ / \sqrt{x(1-x)}$.

**6. Discussion and some further problems.** In Section 4 we have shown that the statistics $\{S_n(s) : s \in \mathbb{R}\}$ behave quite differently for $0 < s < 1$, for $s \le 0$ and $s \ge 1$. In particular, for $0 < s < 1$, the statistics $S_n(s)$ converge almost surely to their natural parameters for fixed $F \ne F_0$. Moreover, the different Poisson boundary behaviors for $s < 0$ and $s \ge 1$ suggest that the statistics $S_n(s)$ with $s \ge 1$ are geared toward "heavy tails," while the statistics $S_n(s)$ with $s \le 0$ are geared more toward "light tails." This also becomes apparent from plots of the functions $K_1(F_1(x), x)$, $K_1(\tilde{F}_0(x), x)$, and of $K_0(F_1(x), x)$ and $K_0(\tilde{F}_0(x), x)$, where $F_1(x) = 1/(1 + \log(1/x))$ and $\tilde{F}_0(x) = \exp(-(1/x - 1))$; see, for example, Jager [23], page 11.

In Section 5 we have shown that all of the statistics $\{S_n(s) : -1 \le s \le 2\}$ have the same optimal detection boundary for the two-point normal mixture testing problem considered by Donoho and Jin [15]. Thus, we have some flexibility in designing a test to detect these subtle tail alternatives, and yet behaving very stably under fixed alternatives (in the sense of always consistently estimating a natural parameter). Thus, it seems that the (Hellinger-type) statistic $S_n(1/2)$ may be a very reasonable compromise test statistic.

Here is a brief listing of some of the remaining open problems:

- Problem 1: What is the limit distribution of $S_n(s)$ under the null hypothesis when $s < -1$ or $s > 2$?

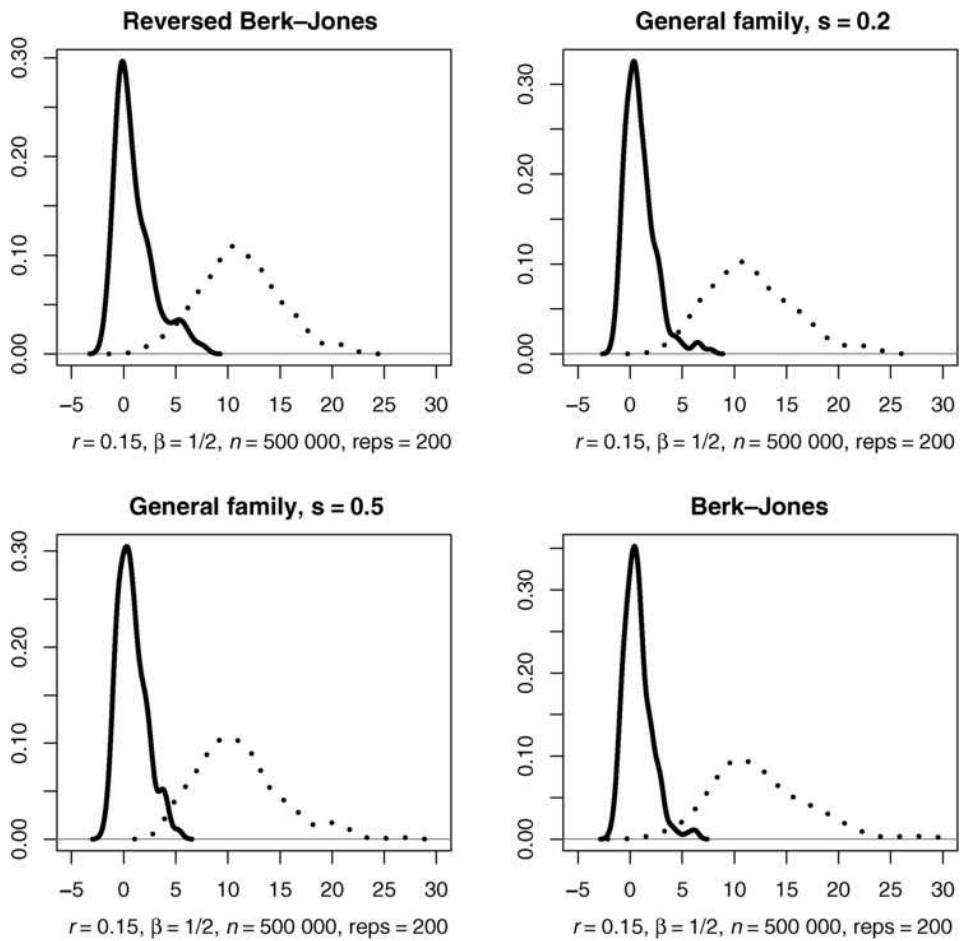

Fig. 1.  *Smoothed histograms of (reps) values of the statistics* $S_n^+(s)$ *under the null hypothesis* $H_0$ *(solid line) and alternative hypothesis* $H_1^{(n)}$ *(dotted line) with* $r = 0.15$, $\beta = 1/2$ *for* $s = 0, 0.2, 0.5, 1$.

- Problem 2: What are necessary and sufficient conditions for $S_n(s)$ to converge to its natural parameter under fixed alternatives $F$ for $s \leq 0$?
- Problem 3: What is the large deviation behavior of $S_n(s)$ under the null hypothesis for $0 < s < 1$?
- Problem 4: Is there an appropriate contiguity theory for the statistics $S_n(s)$? (The only example involving something similar of which we are aware is Theorem A1 of Bickel and Rosenblatt [6], but their results do not seem to apply to the statistics $S_n(s)$.)
- It is fairly easy to construct versions of our statistics $S_n(s)$ in more general settings by replacing the intervals $[0, x]$ and $(x, 1]$ with sets $C$ and $C^c$ for $C$ in some class of sets $\mathcal{C}$. Then for testing $H_0 : P = P_0$ versus $H_1 : P_1 \neq P_0$, a natural generalization of the statistics $S_n(s)$ is

$$S_n(s, \mathcal{C}) = \sup_{C \in \mathcal{C}} K_s(\mathbb{P}_n(C), P_0(C)),$$

where $\mathbb{P}_n$ is the empirical measure of $X_1, \dots, X_n$ i.i.d. $P$.

- Problem 5: Do the statistics $S_n(s, \mathcal{C})$ have reasonable power behavior for some of the "chimeric alternatives" of Khmaladze [30] for some choice of $\mathcal{C}$?

**7. Proofs.**

7.1. *Proofs for Section 3.*

PROOF OF THEOREM 3.1.  We first carry out the proof for $-1 \leq s < 1$, and then indicate the changes that are necessary for $1 \leq s \leq 2$. Fix $s \in [-1, 1)$. Note that

$$\frac{\partial}{\partial u} K_s(u, v) \Big|_{u=v} = \phi_s\left(\frac{u}{v}\right) - \phi_s\left(\frac{1-u}{1-v}\right)\Big|_{u=v} = \phi_s(1) - \phi_s(1) = 0$$

and

$$\frac{\partial^2}{\partial u^2} K(u, v) = \left(\frac{u}{v}\right)^{s-2} \frac{1}{v} + \left(\frac{1-u}{1-v}\right)^{s-2} \frac{1}{1-v} \equiv D_s(u, v).$$

Hence, it follows by Taylor expansion of $u \mapsto K_s(u, v)$ about $u = v$ that

$$K_s(u, v) = K_s(v, v) + \frac{\partial}{\partial u} K_s(u, v) \Big|_{u=v} (u - v) + \frac{1}{2} \frac{\partial^2}{\partial u^2} K_s(u, v) \Big|_{u=u^*} (u - v)^2$$

$$= 0 + 0 + \frac{1}{2}(u - v)^2 D_s(u^*, v) = \frac{1}{2}(u - v)^2 D_s(u^*, v)$$

for some $u^*$ satisfying $|u^* - v| \leq |u - v|$. This yields

$$(9) \qquad K_s(\mathbb{F}_n(x), x) = \tfrac{1}{2}(\mathbb{F}_n(x) - x)^2 D_s(\mathbb{F}_n^*(x), x)$$

for $0 < x < 1$, where $|\mathbb{F}_n^*(x) - x| \leq |\mathbb{F}_n(x) - x|$; that is, $x \leq \mathbb{F}_n^*(x) \leq \mathbb{F}_n(x)$ on the event $x \leq \mathbb{F}_n(x)$ and $\mathbb{F}_n(x) \leq \mathbb{F}_n^*(x) \leq x$ on the event $\mathbb{F}_n(x) \leq x$.

We can write (9) as

$$(10) \qquad K_s(\mathbb{F}_n(x), x) = \frac{1}{2} \frac{(\mathbb{F}_n(x) - x)^2}{x(1 - x)} \{1 + x(1 - x) D_s(\mathbb{F}_n^*(x), x) - 1\},$$

where

$$|\text{Rem}_n(x)| \equiv |x(1 - x) D_s(\mathbb{F}_n^*(x), x) - 1|$$

$$= \left| x(1 - x) \left\{ \left(\frac{\mathbb{F}_n^*(x)}{x}\right)^{s-2} \frac{1}{x} + \left(\frac{1 - \mathbb{F}_n^*(x)}{1 - x}\right)^{s-2} \frac{1}{1 - x} \right\} - 1 \right|$$

$$= \left| (1-x)\left(\frac{x}{\mathbb{F}_n^*(x)}\right)^{2-s} + x\left(\frac{1-x}{1-\mathbb{F}_n^*(x)}\right)^{2-s} - (1-x) - x \right|$$

$$\leq (1-x)\left|\left(\frac{x}{\mathbb{F}_n^*(x)}\right)^{2-s} - 1\right| + x\left|\left(\frac{1-x}{1-\mathbb{F}_n^*(x)}\right)^{2-s} - 1\right|.$$

Fix $\delta \in (0, 1/2)$. Now for $x \in [\delta, 1-\delta]$, $\mathbb{F}_n(x) \in [\delta/2, 1-\delta/2]$ a.s. for $n \geq N_\omega$, so, much as in Wellner and Koltchinskii [43],

$$\sup_{\delta \leq x \leq 1-\delta} |\mathrm{Rem}_n(x)| = O_p(n^{-1/2}).$$

For $0 < v \leq 1/2$ and $1 \leq s \leq 2$, the function $u \mapsto D_s(u,v)$ is monotone for $u \in (0, 1/2]$, while for $0 < v \leq 1/2$ and $-1 \leq s \leq 1$, the function $u \mapsto D_s(u,v)$ is monotone for $u \in (0, b(v,s)]$, where $b(v,s) \equiv 1/(1+c(v,s))$ with $c(v,s) \equiv ((1-v)/v)^{(1-s)/(3-s)}$, so

$$x(1-x)D_s(\mathbb{F}_n^*(x), x) \leq x(1-x)\{D_s(x,x) \vee D_s(\mathbb{F}_n(x), x)\}$$

on the set $\{\mathbb{F}_n(x) < 1/2 \wedge b(x,s)\}$. Since $P(\mathbb{F}_n(\delta) \geq 1/2 \wedge b(\delta,s)) \to 0$ for $\delta < 1/2$ and $s \in [-1,2]$, we get

$$\begin{aligned}
\sup_{X_{(1)} \leq x \leq \delta} |\mathrm{Rem}_n(x)| &\leq \sup_{X_{(1)} \leq x \leq \delta} \left|\left(\frac{x}{\mathbb{F}_n(x)}\right)^{2-s} - 1\right| \\
&\quad + \sup_{X_{(1)} \leq x \leq \delta} \left|\left(\frac{1-x}{1-\mathbb{F}_n(x)}\right)^{2-s} - 1\right| \\
&= O_p(1) + o_p(1) = O_p(1)
\end{aligned} \tag{11}$$

by Shorack and Wellner [38], inequality 1, page 415, and inequality 2, (10.3.6), page 416. Alternatively, see Wellner [42], Lemma 2, page 75, and Remark 1(ii). Similarly,

$$\begin{aligned}
\sup_{1-\delta \leq x < X_{(n)}} |\mathrm{Rem}_n(x)| &\leq \sup_{1-\delta \leq x < X_{(n)}} \left|\left(\frac{x}{\mathbb{F}_n(x)}\right)^{2-s} - 1\right| \\
&\quad + \sup_{1-\delta \leq x < X_{(n)}} \left|\left(\frac{1-x}{1-\mathbb{F}_n(x)}\right)^{2-s} - 1\right| \\
&= o_p(1) + O_p(1) = O_p(1).
\end{aligned} \tag{12}$$

Now we write

$$S_n(s) = S_n(s, I) \vee S_n(s, II) \vee S_n(s, III),$$

where

$$S_n(s, I) \equiv \sup_{\delta \leq x \leq 1-\delta} K_s(\mathbb{F}_n(x), x) = \sup_{\delta \leq x \leq 1-\delta} \tfrac{1}{2}(\mathbb{F}_n(x) - x)^2 D_s(\mathbb{F}_n^*(x), x),$$

$$S_n(s, II) \equiv \sup_{X_{(1)} \leq x \leq \delta} K_s(\mathbb{F}_n(x), x) = \tfrac{1}{2} \sup_{X_{(1)} \leq x \leq \delta} (\mathbb{F}_n(x) - x)^2 D_s(\mathbb{F}_n^*(x), x),$$

$$S_n(s, III) \equiv \sup_{1-\delta \leq x < X_{(n)}} K_s(\mathbb{F}_n(x), x) = \tfrac{1}{2} \sup_{1-\delta \leq x < X_{(n)}} (\mathbb{F}_n(x) - x)^2 D_s(\mathbb{F}_n^*(x), x).$$

By the monotonicity of $u \mapsto D_s(u, v)$ for $u \leq 1/2$ again, with probability tending to 1,

$$S_n(s, II) \leq \tfrac{1}{2} \sup_{X_{(1)} \leq x \leq \delta} (\mathbb{F}_n(x) - x)^2 \{ D_s(\mathbb{F}_n(x), x) \vee D_s(x, x) \}$$

$$\geq \tfrac{1}{2} \sup_{X_{(1)} \leq x \leq \delta} (\mathbb{F}_n(x) - x)^2 \{ D_s(\mathbb{F}_n(x), x) \wedge D_s(x, x) \},$$

and similarly, by the monotonicity of $u \mapsto D_s(u, v)$ for $1/2 \leq u < 1$,

$$S_n(s, III) \leq \tfrac{1}{2} \sup_{1-\delta \leq x < X_{(n)}} (\mathbb{F}_n(x) - x)^2 \{ D_s(\mathbb{F}_n(x), x) \vee D_s(x, x) \}$$

$$\geq \tfrac{1}{2} \sup_{1-\delta \leq x < X_{(n)}} (\mathbb{F}_n(x) - x)^2 \{ D_s(\mathbb{F}_n(x), x) \wedge D_s(x, x) \}.$$

For $S_n(s, I)$,

$$S_n(s, I) = \frac{1}{2} \sup_{\delta \leq x \leq 1-\delta} \frac{(\mathbb{F}_n(x) - x)^2}{x(1-x)} \{ 1 + O_p(n^{-1/2}) \}$$

$$\leq \frac{1}{2} \sup_{\delta \leq x \leq 1-\delta} (\mathbb{F}_n(x) - x)^2 \{ D_s(\mathbb{F}_n(x), x) \vee D_s(x, x) \} \{ 1 + O_p(n^{-1/2}) \}.$$

In the second region the argument above leading to (11) yields

$$S_n(s, II) \leq \tfrac{1}{2} \sup_{X_{(1)} \leq x \leq \delta} (\mathbb{F}_n(x) - x)^2 \{ D_s(\mathbb{F}_n(x), x) \vee D_s(x, x) \}$$

$$\geq \tfrac{1}{2} \sup_{X_{(1)} \leq x \leq \delta} (\mathbb{F}_n(x) - x)^2 \{ D_s(\mathbb{F}_n(x), x) \wedge D_s(x, x) \},$$

and similarly for $S_n(s, III)$. It follows that

$$S_n(s) \leq \frac{1}{2} \sup_{X_{(1)} \leq x < X_{(n)}} (\mathbb{F}_n(x) - x)^2$$

$$\times \{ D_s(\mathbb{F}_n(x), x) \vee D_s(x, x) \} \{ 1 + O_p(n^{-1/2}) \}$$

(13)

$$= \frac{1}{2} \sup_{X_{(1)} \leq x < X_{(n)}} \left\{ \frac{(\mathbb{F}_n(x) - x)^2}{x(1-x)} \{ 1 \vee x(1-x) D_s(\mathbb{F}_n(x), x) \} \right\}$$

$$\times \{ 1 + O_p(n^{-1/2}) \}$$

and, on the other hand,

$$\begin{aligned}
S_n(s) &\geq \frac{1}{2} \sup_{X_{(1)} \leq x < X_{(n)}} (\mathbb{F}_n(x) - x)^2 \{D_s(\mathbb{F}_n(x), x) \wedge D_s(x, x)\} \\
&\quad \times \{1 + O_p(n^{-1/2})\} \\
&= \frac{1}{2} \sup_{X_{(1)} \leq x < X_{(n)}} \left\{ \frac{(\mathbb{F}_n(x) - x)^2}{x(1-x)} \{1 \wedge x(1-x) D_s(\mathbb{F}_n(x), x)\} \right\} \\
&\quad \times \{1 + O_p(n^{-1/2})\}.
\end{aligned}$$

(14)

Now we break the supremum into the regions $[X_{(1)}, d_n]$, $[d_n, 1 - d_n]$ and $[1 - d_n, X_{(n)})$ with $d_n = (\log n)^k / n$ for any $k \geq 1$. Then we have

$$n \sup_{X_{(1)} \leq x \leq d_n} \frac{(\mathbb{F}_n(x) - x)^2}{x(1-x)} = o_p(b_n^2),$$

where $b_n = \sqrt{2 \log_2 n}$; see Shorack and Wellner [38], (26), page 602. Moreover,

$$\sup_{X_{(1)} \leq x \leq d_n} |x(1-x) D_s(\mathbb{F}_n(x), x)| = O_p(1),$$

so

(15)    $$n \sup_{X_{(1)} \leq x \leq d_n} \frac{(\mathbb{F}_n(x) - x)^2}{x(1-x)} (1 \# x(1-x) D_s(\mathbb{F}_n(x), x)) = o_p(b_n^2)$$

for $\# = \wedge$ or $\# = \vee$, and similarly for the region $[1 - d_n, X_{(n)}]$ by a symmetric argument. On the other hand, if we define

$$Z_n \equiv \sup_{d_n \leq x \leq 1 - d_n} \frac{\sqrt{n}|\mathbb{F}_n(x) - x|}{\sqrt{x(1-x)}},$$

then, for $k \geq 5$,

(16)                        $$\frac{Z_n}{b_n} \xrightarrow{p} 1$$

and

(17)                        $$b_n Z_n - c_n \xrightarrow{d} Y_4 \sim E_v^4,$$

where $c_n = 2 \log_2 n + (1/2) \log_3 n - (1/2) \log(4\pi)$ (see, e.g., Shorack and Wellner [38], page 600, (16.1.20) and (16.1.17)). [Note that for the middle bracket in (15) we have

$$x(1-x) D_s(\mathbb{F}_n(x), x) = (1-x)\left(\frac{x}{\mathbb{F}_n(x)}\right)^{2-s} + x\left(\frac{1-x}{1 - \mathbb{F}_n(x)}\right)^{2-s},$$

so

$$\sup_{X_{(1)} \leq x \leq X_{(1)}} x(1-x) D_s(\mathbb{F}_n(x), x) \leq \left( \sup_{X_{(1)} \leq x \leq d_n} \frac{x}{\mathbb{F}_n(x)} \right)^{2-s} + o_p(1),$$

so by Wellner [42], Remark 1, the probability of large values of the main term can be bounded by

$$P\left( \left( \sup_{X_{(1)} \leq x \leq d_n} \frac{x}{\mathbb{F}_n(x)} \right)^{2-s} > \lambda \right) = P\left( \sup_{X_{(1)} \leq x \leq d_n} \frac{x}{\mathbb{F}_n(x)} > \lambda^{1/(2-s)} \right)$$

$$\leq P\left( \sup_{X_{(1)} \leq x \leq 1} \frac{x}{\mathbb{F}_n(x)} > \lambda^{1/(2-s)} \right)$$

$$\leq e \cdot \lambda^{1/(2-s)} \exp(-\lambda^{1/(2-s)}). \Big]$$

Furthermore,

$$(18) \qquad \left\| \frac{\mathbb{F}_n(x) - x}{x} \right\|_{d_n}^1 = O(a_n)$$

almost surely where

$$a_n^2 \equiv \frac{\log_2 n}{n d_n} = \frac{\log_2 n}{(\log n)^k} \to 0;$$

see Shorack and Wellner [38], page 424, (4.5.10) and (4.5.11). It follows from (13), (14) and (15)–(18) that

$$nS_n(s) = \frac{1}{2} \left\{ \sup_{d_n \leq x \leq 1 - d_n} \frac{n(\mathbb{F}_n(x) - x)^2}{x(1-x)} (1 + O_p(a_n)) \vee o_p(b_n^2) \right\}$$

$$(19) \qquad \times \{1 + O_p(n^{-1/2})\}$$

$$= \frac{1}{2} \{ Z_n^2 \vee o_p(b_n^2) \} + o_p(1).$$

Hence, we can write

$$\frac{1}{2} Z_n^2 = \frac{1}{2} (Z_n - c_n/b_n)(Z_n + c_n/b_n) + \frac{1}{2} \frac{c_n^2}{b_n^2}$$

$$= \frac{1}{2} b_n (Z_n - c_n/b_n) \frac{Z_n + c_n/b_n}{b_n} + \frac{1}{2} \frac{c_n^2}{b_n^2}.$$

It follows that

$$nS_n(s) - \frac{1}{2} \frac{c_n^2}{b_n^2}$$

$$(20) \quad = b_n(Z_n - c_n/b_n)\frac{Z_n + c_n/b_n}{2b_n} \vee \left(o_p(b_n^2) - \frac{1}{2}\frac{c_n^2}{b_n^2}\right) + o_p(1)$$

$$= b_n(Z_n - c_n/b_n)\frac{Z_n/b_n + c_n/b_n^2}{2} \vee (o_p(1) - 1/2)b_n^2 + o_p(1)$$

$$\xrightarrow{d} Y_4\frac{1+1}{2} \vee \{-\infty\} = Y_4;$$

here we used $c_n^2/b_n^2 \sim b_n^2$ in the second equality. Since

$$\frac{1}{2}\frac{c_n^2}{b_n^2} = \log_2 n + (1/2)\log_3 n - (1/2)\log(4\pi) + o(1) = r_n + o(1),$$

this yields

$$(21) \quad P(nS_n(s) - r_n \leq x) \rightarrow \exp(-4\exp(-x)),$$

and completes the proof of Theorem 3.1. Note that the centering $c_n^2/(2b_n^2)$ emerges naturally in the course of this proof. This completes the proof for the case $s \in [-1, 1)$. For $1 \leq s \leq 2$, there are two additional terms that enter, and both of these are $o_p(b_n^2)$ from the arguments in the previous section. The case $s = 2$ is easy since in this case $v(1-v)D_s(u, v) = 1$ for all $u$, while the result was stated for the case $s = 1$ by Berk and Jones [4] and proved in Wellner and Koltchinskii [43]. (Wellner and Koltchinskii [43] incorrectly claim (page 324) that $K(\mathbb{F}_n(x), x) = 0$ if $x < X_{(1)}$; in fact, the supremum over this region is stochastically bounded and, hence, can be neglected.) □

PROOF OF THEOREM 3.2. The fact that $nT_n(2) = A_n^2/2 \xrightarrow{d} A^2/2$ is classical; see Shorack and Wellner [38], page 148. That $nT_n(1) \xrightarrow{d} A^2/2$ was noted by Einmahl and McKeague [18] and proved by Wellner and Koltchinskii [43]. The proof for $s \neq 1, 2$ proceeds along the same lines as the proof in Wellner and Koltchinskii [43] for the case $s = 1$, and hence will not be given here. For details, see Jager [23] or Jager and Wellner [26]. □

7.2. *Proofs for Section 4.*

PROOF OF PROPOSITION 4.1. We first prove the claim for the "unrestricted version" of the statistics $S_n^{ur}(s)$ defined by $S_n^{ur}(s) \equiv \sup_{0 < x < 1} K_s \times (\mathbb{F}_n(x), x)$, and then show that the difference between $S_n(s)$ and $S_n^{ur}(s)$ is negligible. Now for $s \in (0, 1)$ and $C_s \equiv 1/(s(1-s))$, we have

$$|S_n^{ur}(s) - S_\infty(s, F)|$$

$$\leq C_s \sup_{0 < x < 1} |\{1 - \mathbb{F}_n(x)^s x^{1-s} - (1 - \mathbb{F}_n(x))^s(1-x)^{1-s}\}$$

$$- \{1 - F(x)^s x^{1-s} - (1 - F(x))^s(1-x)^{1-s}\}|$$

$$\leq C_s \left\{ \sup_x |(\mathbb{F}_n(x)^s - F(x)^s) x^{1-s}| \right.$$

$$+ \sup_x |\{(1 - \mathbb{F}_n(x))^s - (1 - F(x))^s\} x^{1-s}| \left. \right\}$$

$$\leq C_s \left\{ \sup_x |\mathbb{F}_n(x)^s - F(x)^s| + \sup_x |(1 - \mathbb{F}_n(x))^s - (1 - F(x))^s| \right\}$$

$$\overset{a.s.}{\to} 0.$$

Thus, the proposition will be proved if we show that

(22) $$S_n^{ur}(s) - S_n(s) \overset{a.s.}{\to} 0.$$

Now write $S_n^0(s) = \max\{R_n, M_n, L_n\}$, where

$$M_n \equiv \sup_{X_{(1)} \leq x < X_{(n)}} K_s(\mathbb{F}_n(x), x) = S_n(s),$$

$$L_n \equiv \sup_{x < X_{(1)}} K_s(\mathbb{F}_n(x), x)$$

and

$$R_n \equiv \sup_{x \geq X_{(n)}} K_s(\mathbb{F}_n(x), x).$$

Note that

$$S_n^0(s) - S_n(s) = \max\{L_n, M_n, R_n\} - M_n = \begin{cases} 0, & \text{if } M_n \geq L_n \vee R_n, \\ L_n - M_n, & \text{if } L_n > M_n \vee R_n, \\ R_n - M_n, & \text{if } R_n > M_n \vee L_n. \end{cases}$$

Now set

$$\alpha_0 \equiv \alpha_0(F) = \sup\{x : F(x) = 0\} \geq 0,$$

$$\alpha_1 \equiv \alpha_1(F) = \inf\{x : F(x) = 1\} \leq 1.$$

Note that

$$L_n = \sup_{x < X_{(1)}} K_s(\mathbb{F}_n(x), x) = \frac{1}{s(1-s)} \{1 - (1 - X_{(1)})^{1-s}\}$$

$$\overset{a.s.}{\to} \frac{1}{s(1-s)} \{1 - (1 - \alpha_0)^{1-s}\} \equiv l_0(s, F),$$

and, on the other hand,

$$M_n = \sup_{X_{(1)} \leq x < X_{(n)}} K_s(\mathbb{F}_n(x), x) \geq K_s(\mathbb{F}_n(X_{(1)}), X_{(1)}) = K_s(1/n, X_{(1)})$$

$$= \frac{1}{s(1-s)} \{1 - (1/n)^s X_{(1)}^{1-s} - (1 - 1/n)^s (1 - X_{(1)})^{1-s}\} \equiv L_n^0$$

$$\overset{a.s.}{\to} \frac{1}{s(1-s)} \{1 - (1 - \alpha_0)^{1-s}\} = l_0(s, F).$$

Similarly,

$$R_n = \sup_{x \geq X_{(n)}} K_s(\mathbb{F}_n(x), x) = \frac{1}{s(1-s)}\{1 - X_{(n)}^{1-s}\}$$

$$\xrightarrow{a.s.} \frac{1}{s(1-s)}\{1 - \alpha_1^{1-s}\} \equiv r_0(s, F),$$

while

$$M_n = \sup_{X_{(1)} \leq x < X_{(n)}} K_s(\mathbb{F}_n(x), x) \geq K_s(\mathbb{F}_n(X_{(n)}-), X_{(n)}) = K_s(1 - 1/n, X_{(n)})$$

$$= \frac{1}{s(1-s)}\{1 - (1 - 1/n)^s X_{(n)}^{1-s} - (1/n)^s(1 - X_{(n)})^{1-s}\} \equiv R_n^0$$

$$\xrightarrow{a.s.} \frac{1}{s(1-s)}\{1 - \alpha_1^{1-s}\} = r_0(s, F).$$

By combining these pieces, it follows that

$$0 \leq S_n^0(s) - S_n(s) = \max\{L_n, M_n, R_n\} - M_n$$

$$= \left\{\begin{array}{ll} 0, & \text{if } M_n \geq L_n \vee R_n, \\ L_n - M_n, & \text{if } L_n > M_n \vee R_n, \\ R_n - M_n, & \text{if } R_n > M_n \vee L_n \end{array}\right\} \leq \left\{\begin{array}{ll} 0, & \text{if } M_n \geq L_n \vee R_n, \\ L_n - L_n^0, & \text{if } L_n > L_n^0 \vee R_n, \\ R_n - R_n^0, & \text{if } R_n > R_n^0 \vee L_n \end{array}\right\}$$

$$\xrightarrow{a.s.} 0.$$

This shows that (22) holds and completes the proof. □

PROOF OF PROPOSITION 4.2. Recall that when $s = 2$ such a condition follows from Theorem 3 in Jager and Wellner [25]: taking $b = 1/2$ and applying continuous mapping, we conclude that

$$S_n(2) \xrightarrow{a.s.} \sup_{0 < x < 1} K_2(F(x), x) \quad \text{if and only if} \quad E\{[X(1 - X)]^{-1/2}\} < \infty.$$

Similarly, for $1 < s < \infty$,

$$S_n(s) = \sup_{0 < x < 1}\{\mathbb{F}_n(x)^s x^{1-s} + (1 - \mathbb{F}_n(x))^s(1 - x)^{1-s} - 1\}\frac{1}{s(s-1)}$$

$$\xrightarrow{a.s.} \sup_{0 < x < 1}\{F(x)^s x^{1-s} + (1 - F(x))^s(1 - x)^{1-s} - 1\}\frac{1}{s(s-1)}$$

if and only if

$$\|\mathbb{F}_n(x)^s x^{1-s} - F(x)^s x^{1-s}\| \xrightarrow{a.s.} 0$$

and

$$\|(1 - \mathbb{F}_n(x))^s(1 - x)^{1-s} - (1 - F(x))^s(1 - x)^{1-s}\| \xrightarrow{a.s.} 0,$$

if and only if

$$\left\|\left(\frac{\mathbb{F}_n(x)}{x^{(s-1)/s}}\right)^s - \left(\frac{F(x)}{x^{(s-1)/s}}\right)^s\right\| \overset{a.s.}{\to} 0$$

and

$$\left\|\left(\frac{1-\mathbb{F}_n(x)}{(1-x)^{(s-1)/s}}\right)^s - \left(\frac{1-F(x)}{(1-x)^{(s-1)/s}}\right)^s\right\| \overset{a.s.}{\to} 0,$$

if and only if

$$\left\|\left(\frac{\mathbb{G}_n(F(x))}{x^{(s-1)/s}}\right)^s - \left(\frac{F(x)}{x^{(s-1)/s}}\right)^s\right\| \overset{a.s.}{\to} 0$$

and

$$\left\|\left(\frac{1-\mathbb{G}_n(F(x))}{(1-x)^{(s-1)/s}}\right)^s - \left(\frac{1-F(x)}{(1-x)^{(s-1)/s}}\right)^s\right\| \overset{a.s.}{\to} 0.$$

Since $g(u) = u^s$ is uniformly continuous on bounded sets, these last two convergences occur if and only

$$\left\|\frac{\mathbb{G}_n(F(x))}{x^{(s-1)/s}} - \frac{F(x)}{x^{(s-1)/s}}\right\| \overset{a.s.}{\to} 0$$

and

$$\left\|\frac{1-\mathbb{G}_n(F(x))}{(1-x)^{(s-1)/s}} - \frac{1-F(x)}{(1-x)^{(s-1)/s}}\right\| \overset{a.s.}{\to} 0.$$

These, in turn, hold if and only if

$$\left\|\frac{\mathbb{G}_n(u)-u}{F^{-1}(u)^{(s-1)/s}}\right\| \overset{a.s.}{\to} 0 \quad\text{and}\quad \left\|\frac{1-\mathbb{G}_n(u)-(1-u)}{(1-F^{-1}(u))^{(s-1)/s}}\right\| \overset{a.s.}{\to} 0.$$

But in view of Wellner [41], these convergences hold if and only if $F$ satisfies

$$\int_0^1 \frac{1}{(F^{-1}(u)(1-F^{-1}(u)))^{(s-1)/s}}\,du < \infty.$$

By the (inverse) probability integral transformation, the convergence in the last display is equivalent to $E[X(1-X)]^{(1-s)/s} < \infty$. This completes the proof of the claimed equivalences.  □

PROOF OF PROPOSITION 4.3.  For $s = 1$, this follows from Berk and Jones [5], pages 55–56. Thus, it suffices to prove the claimed convergences for $s > 1$ and $s < 0$.

For $s > 1$, fix $\alpha \in (0, 1)$. We begin by breaking the supremum over $(0, 1)$ into the regions $0 < F_s(x) < n^{-\alpha}$, $n^{-\alpha} \le F_s(x) \le 1 - n^{-\alpha}$ and $1 - n^{-\alpha} < F_s(x) < 1$:

$$S_n(s) = \sup_{0 < x < 1} \{\mathbb{F}_n(x)^s x^{1-s} + (1 - \mathbb{F}_n(x))^x (1-x)^{1-s} - 1\} \frac{1}{s(s-1)}$$

$$= \sup_{x:0 < F_s(x) < n^{-\alpha}} \{\mathbb{F}_n(x)^s x^{1-s} + (1 - \mathbb{F}_n(x))^x (1-x)^{1-s} - 1\} \frac{1}{s(s-1)}$$

$$\vee \sup_{x:n^{-\alpha} < F_s(x) < 1 - n^{-\alpha}} \{\mathbb{F}_n(x)^s x^{1-s} + (1 - \mathbb{F}_n(x))^x (1-x)^{1-s} - 1\} \frac{1}{s(s-1)}$$

$$\vee \sup_{x:1 - n^{-\alpha} < F_s(x) < 1} \{\mathbb{F}_n(x)^s x^{1-s} + (1 - \mathbb{F}_n(x))^x (1-x)^{1-s} - 1\} \frac{1}{s(s-1)}$$

$$\equiv I_n(s) \vee II_n(s) \vee III_n(s).$$

For the main term, $I_n(s)$, let $\mathbb{G}_n$ be the empirical d.f. of $n$ i.i.d. Uniform$(0,1)$ random variables and use $\mathbb{F}_n \stackrel{d}{=} \mathbb{G}_n(F_s)$ to write

$$s(s-1)I_n(s) \stackrel{d}{=} \sup_{0 < F_s(x) < n^{-\alpha}} \left\{ \left( \frac{\mathbb{G}_n(F_s(x))}{F_s(x)} \right)^s F_s(x)^s x^{1-s} \right.$$

$$\left. + \left( \frac{1 - \mathbb{G}_n(F_s(x))}{1 - F_s(x)} \right)^s (1 - F_s(x))^s (1-x)^{1-s} - 1 \right\},$$

where

$$\sup_{x:F_s(x) < n^{-\alpha}} \left( \frac{\mathbb{G}_n(F_s(x))}{F_s(x)} \right)^s = \sup_{0 < t < n^{-\alpha}} \left( \frac{\mathbb{G}_n(t)}{t} \right)^s \stackrel{d}{\to} \sup_{t>0} \left( \frac{\mathbb{N}(t)}{t} \right)^s,$$

$$F_s(x)^s x^{1-s} = \frac{x^{1-s}}{1 + (x^{1-s} - 1)/(s-1)} \to s - 1$$

uniformly in $x \in [0, n^{-\alpha}]$, while

$$\sup_{x:F_s(x) < n^{-\alpha}} \left| \frac{1 - \mathbb{G}_n(F_s(x))}{1 - F_s(x)} - 1 \right| \stackrel{a.s.}{\to} 0$$

and

$$\sup_{x:F_s(x) < n^{-\alpha}} |(1 - F_s(x))^s (1-x)^{1-s} - 1| \to 0.$$

Combining these last five displays shows that $I_n(s) \stackrel{d}{\to} s^{-1} \sup_{t>0}(\mathbb{N}(t)/t)^s$; note that the limit variable is $\ge 1/s$ almost surely.

To handle the term $II_n(s)$, write

$$II_n(s) \overset{d}{=} \sup_{n^{-\alpha} < F_s(x) < 1-n^{-\alpha}} \left\{ \left( \frac{\mathbb{G}_n(F_s(x))}{F_s(x)} \right)^s F_s(x)^s x^{1-s} \right.$$
$$+ \left( \frac{1 - \mathbb{G}_n(F_s(x))}{1 - F_s(x)} \right)^s$$
$$\left. \times (1 - F_s(x))^s (1-x)^{1-s} - 1 \right\} \frac{1}{s(s-1)},$$

where now the two terms involving the ratio of the empirical d.f. to the true d.f. $F_s$ converge almost surely to 1. Hence, we conclude that

$$II_n(s) \overset{a.s.}{\longrightarrow} \sup_{0<x<1} K_s(F_s(x),x) = \frac{1}{s},$$

where the equality follows after some calculation. Finally, it is easily shown that $III_n(s) \overset{a.s.}{\longrightarrow} 0$.

For $s < 0$, fix $\alpha \in (0,1)$. We begin by breaking the supremum over $(0,1)$ into the regions $X_{(1)} \le x < F_n^{-1}(n^{-\alpha})$, $F_s^{-1}(n^{-\alpha}) \le x \le F_s^{-1}(1-n^{-\alpha})$ and $F_s^{-1}(1-n^{-\alpha}) < x < X_{(n)}$:

$$S_n(s) = \sup_{X_{(1)} \le x < X_{(n)}} \{ \mathbb{F}_n(x)^s x^{1-s} + (1 - \mathbb{F}_n(x))^x (1-x)^{1-s} - 1 \} \frac{1}{s(s-1)}$$

$$= \sup_{x:X_{(1)} \le x < F_s^{-1}(n^{-\alpha})} \{ \mathbb{F}_n(x)^s x^{1-s} + (1 - \mathbb{F}_n(x))^x (1-x)^{1-s} - 1 \} \frac{1}{s(s-1)}$$

$$\vee \sup_{x:F_s^{-1}(n^{-\alpha}) \le x \le F_s^{-1}(1-n^{-\alpha})} \{ \mathbb{F}_n(x)^s x^{1-s} + (1 - \mathbb{F}_n(x))^x (1-x)^{1-s} - 1 \} \frac{1}{s(s-1)}$$

$$\vee \sup_{x:F_s^{-1}(1-n^{-\alpha}) < x < X_{(n)}} \{ \mathbb{F}_n(x)^s x^{1-s} + (1 - \mathbb{F}_n(x))^x (1-x)^{1-s} - 1 \} \frac{1}{s(s-1)}$$

$$\equiv I_n(s) \vee II_n(s) \vee III_n(s).$$

For the main term, $I_n(s)$, let $\mathbb{G}_n$ be the empirical d.f. of $n$ i.i.d. Uniform$(0,1)$ random variables and use $\mathbb{F}_n \overset{d}{=} \mathbb{G}_n(F_s)$ to write

$$s(s-1)I_n(s) \overset{d}{=} \sup_{x:X_{(1)} \le x < F_s^{-1}(n^{-\alpha})} \left\{ \left( \frac{F_s(x)}{\mathbb{G}_n(F_s(x))} \right)^{-s} F_s(x)^s x^{1-s} \right.$$

$$+ \left( \frac{1 - F_s(x)}{1 - \mathbb{G}_n(F_s(x))} \right)^{-s}$$

$$\times (1 - F_s(x))^s (1 - x)^{1-s} - 1 \Big\},$$

where

$$\sup_{x : X_{(1)} \leq x < F_s^{-1}(n^{-\alpha})} \left( \frac{F_s(x)}{\mathbb{G}_n(F_s(x))} \right)^{-s} = \sup_{\xi_{(1)} \leq t < n^{-\alpha}} \left( \frac{t}{\mathbb{G}_n(t)} \right)^{-s} \overset{d}{\to} \sup_{t \geq S_1} \left( \frac{t}{\mathbb{N}(t)} \right)^{-s},$$

$$F_s(x)^s x^{1-s} = x^{1-s}(1 - s(x^{-(1-s)} - 1)) \to -s$$

uniformly in $x \in [0, F^{-1}(n^{-\alpha})]$, while

$$\sup_{x : F_s(x) < n^{-\alpha}} \left| \frac{1 - F_s(x)}{1 - \mathbb{G}_n(F_s(x))} - 1 \right| \overset{a.s.}{\to} 0$$

and

$$\sup_{x : F_s(x) < n^{-\alpha}} |(1 - F_s(x))^s (1 - x)^{1-s} - 1| \to 0.$$

Combining these last five displays shows that $I_n(s) \overset{d}{\to} (1-s)^{-1} \sup_{t \geq S_1} (t/\mathbb{N}(t))^{-s}$; note that the limit variable is $\geq 1/(1-s)$ almost surely.

To handle the term $II_n(s)$, write

$$II_n(s) \overset{d}{=} \sup_{n^{-\alpha} < F_s(x) < 1 - n^{-\alpha}} \left\{ \left( \frac{F_s(x)}{\mathbb{G}_n(F_s(x))} \right)^{-s} F_s(x)^s x^{1-s} \right.$$

$$+ \left( \frac{1 - F_s(x)}{1 - \mathbb{G}_n(F_s(x))} \right)^{-s}$$

$$\left. \times (1 - F_s(x))^s (1 - x)^{1-s} - 1 \right\} \frac{1}{s(s-1)},$$

where now the two terms involving the ratio of the empirical d.f. to the true d.f. $F_s$ converge almost surely to 1. Hence, we conclude that

$$II_n(s) \overset{a.s.}{\to} \sup_{0 < x < 1} K_s(F_s(x), x) = \frac{1}{1-s},$$

where the equality follows after some calculation. Finally, it is easily shown that $III_n(s) \overset{a.s.}{\to} 0$. $\quad \square$

To prove Theorem 4.4 and its corollary, we will use the following lemma from Chernoff [9].

LEMMA 7.1. *Let $X_1, X_2, \ldots$ be i.i.d. with continuous distribution $F_0$.*

(a) *If* $t < F_0(x)$, *then* $n^{-1} \log P(\mathbb{F}_n(x) \leq t) \to -K^-(t, F_0(x))$.
(b) *If* $t > F_0(x)$, *then* $n^{-1} \log P(\mathbb{F}_n(x) \geq t) \to -K^+(t, F_0(x))$.

*In both cases, the convergence is from below.*

PROOF.   This follows from Theorem 1 of Chernoff [9].   □

PROOF OF THEOREM 4.4.   We first prove the theorem for $S_n^{ur,+}(s)$. Since $K_s^+(t, x)$ is continuous in $t$ and strictly increasing on $(x, 1)$, then for $0 < a < (1 - x^{1-s})/[s(1-s)]$, there is a unique $\tau = \tau(x)$ in $(x, 1)$ for which $K_s^+(\tau, x) = a$ and $\{t : K_s^+(t, x) \geq a\} = [\tau, \infty)$. If $a \geq (1 - x^{1-s})/[s(1-s)]$, then $\tau = 1$ necessarily.

For any fixed $x \in (0, 1)$, we have

$$\frac{1}{n} \log P(S_n^{ur,+}(s) \geq a) \geq \frac{1}{n} \log P(K_s^+(\mathbb{F}_n(x), x) \geq a)$$

$$= \frac{1}{n} \log P(\mathbb{F}_n(x) \geq \tau) \to -K^+(\tau, x)$$

by Lemma 7.1. So

$$(23) \qquad \liminf_{n \to \infty} \frac{1}{n} \log P(S_n^{ur,+}(s) \geq a) \geq -\inf_{0 < x < 1} K^+(\tau(x), x).$$

Now consider the reverse inequality. Let $\sup_{x|i}$ denote the supremum for $X_{(i)} \leq x < X_{(i+1)}$. Since $\mathbb{F}_n(x) = i/n$ on this range, we have

$$\sup_{x|i} K_s^+(\mathbb{F}_n(x), x) = K_s^+(i/n, X_{(i)}) \vee K_s^+(i/n, X_{(i+1)}) = K_s^+(i/n, X_{(i)}).$$

Note that for $x < X_{(1)}$, we have $\mathbb{F}_n(x) = 0$, and so $K_s^+(\mathbb{F}_n(x), x) = 0$ also. So we can write $S_n^{ur,+}(s) = \max_{1 \leq i \leq n}\{K_s^+(i/n, X_{(i)})\} = \max_{1 \leq i \leq n}\{K_s^+(\mathbb{F}_n(X_{(i)}), X_{(i)})\}$. Now, using monotonicity of $\tau$,

$$\frac{1}{n} \log P(S_n^{ur,+}(s) \geq a)$$

$$= \frac{1}{n} \log P\left(\max_{1 \leq i \leq n}\{K_s^+(\mathbb{F}_n(X_{(i)}), X_{(i)})\} \geq a\right)$$

$$\leq \frac{1}{n} \log \sum_{i=1}^n P(K_s^+(\mathbb{F}_n(X_{(i)}), X_{(i)}) \geq a)$$

$$= \frac{1}{n} \log \sum_{i=1}^n P(\mathbb{F}_n(X_{(i)}) \geq \tau(X_{(i)}))$$

$$= \frac{1}{n} \log \sum_{i=1}^n P(i/n \geq \tau(X_{(i)})) = \frac{1}{n} \log \sum_{i=1}^n P(\tau^{-1}(i/n) \geq X_{(i)})$$

$$\leq \frac{1}{n} \log \sum_{i=1}^{n} P(\mathbb{F}_n(\tau^{-1}(i/n)) \geq \mathbb{F}_n(X_{(i)}))$$

$$= \frac{1}{n} \log \sum_{i=1}^{n} P(\mathbb{F}_n(\tau^{-1}(i/n)) \geq i/n)$$

$$\leq \frac{1}{n} \log \sum_{i=1}^{n} e^{-nK^+(i/n, \tau^{-1}(i/n))} \qquad \text{[by Lemma 7.1(b)]}$$

$$\leq \frac{1}{n} \log \sum_{i=1}^{n} e^{-n\min_{1 \leq i \leq n} K^+(i/n, \tau^{-1}(i/n))} \leq \frac{1}{n} \log \sum_{i=1}^{n} e^{-n\inf_{0<x<1} K^+(x, \tau^{-1}(x))}$$

$$= \frac{1}{n} \log \sum_{i=1}^{n} e^{-n\inf_{0<x<1} K^+(\tau(x), x)} = \frac{1}{n} \log[n e^{-n\inf_{0<x<1} K^+(\tau(x), x)}]$$

$$= -\inf_{0<x<1} K^+(\tau(x), x) + \frac{\log n}{n}.$$

So we conclude that

$$(24) \qquad \limsup_{n \to \infty} \frac{1}{n} \log P(S_n^{ur,+}(s) \geq a) \leq -\inf_{0<x<1} K^+(\tau(x), x).$$

Combining this last display with (23) yields the convergence part of (7). To prove the explicit formula for $g_s^+$, note that $K^+(\tau^+(x), x)$ is a decreasing function of $x$ until $x = [1 - s(1-s)a]^{1/(1-s)}$, where $K^+(\tau^+(x), x) = \infty$. Thus,

$$\inf_{0<x<1} K^+(\tau^+(x), x)$$
$$= \lim_{\epsilon \searrow 0} K^+(\tau^+([1 - s(1-s)a]^{1/(1-s)} - \epsilon), [1 - s(1-s)a]^{1/(1-s)} - \epsilon)$$
$$= -\log[1 - s(1-s)a]/(1-s),$$

so the given formula for the infimum in (7) holds. This completes the proof for $S_n^{ur,+}(s)$. The proof for $S_n^{ur,-}(s)$ is analogous using Lemma 7.1(a).  □

7.3. *Proofs for Section 5.* The following lemma extends Lemma A.4, page 988, Donoho and Jin [15].

LEMMA 7.2.  (i) *For* $0 < v \leq u \leq 1/2$, *and* $-1 \leq s \leq 2$,

$$(25) \qquad K_s^+(u, v) \leq \frac{1}{2} \frac{(u-v)^2}{v(1-v)} \equiv K_2(u, v).$$

(ii) *Let* $1 < s \leq 2$ *and* $v = v(u)$ *satisfy* $0 < v \leq u < 1$. *Then, as* $u \to 0$,

$$K_s(u, v) = \begin{cases} K_2(u, v)[1 + O(1 - (v/u)^{2-s}) \vee O(((1-v)/(1-u))^{2-s} - 1)], \\ \qquad\qquad\qquad\qquad\qquad\qquad\qquad\qquad\qquad \text{if } u/v \to 1, \\ v\phi_s(u/v)(1 + o(1)) = u\{(u/v)^{s-1} - s\}(1 + o(1))/(s(s-1)), \\ \qquad\qquad\qquad\qquad\qquad\qquad\qquad\qquad\qquad \text{if } u/v \to \infty. \end{cases}$$

(iii) *Let* $s = 1$ *and* $v = v(u)$ *satisfy* $0 < v \leq u < 1$. *Then, as* $u \to 0$,

$$K_1(u,v) = \begin{cases} K_2(u,v)[1 + O(u + (u/v) - 1)], & \text{if } u/v \to 1, \\ u \log(u/v)(1 + o(1)), & \text{if } u/v \to \infty. \end{cases}$$

(iv) *Let* $s \in [-1, 1) \setminus \{0\}$ *and* $v = v(u)$ *satisfy* $0 < v \leq u < 1$. *Then, as* $u \to 0$,

$$K_s(u,v) = \begin{cases} K_2(u,v)[1 + O(1 - (v/u)^{2-s}) \vee O(((1-v)/(1-u))^{2-s} - 1)], \\ \qquad \text{if } u/v \to 1, \\ \dfrac{1}{1-s} u(1 + o(1)), \qquad \text{if } u/v \to \infty. \end{cases}$$

(v) *Let* $s = 0$ *and* $v = v(u)$ *satisfy* $0 < v \leq u < 1$. *Then, as* $u \to 0$,

$$K_0(u,v) = \begin{cases} K_2(u,v)[1 + O(1 - (v/u)^2) \vee O(((1-v)/(1-u))^2 - 1)], \\ \qquad \text{if } u/v \to 1, \\ u(1 + o(1)), \qquad \text{if } u/v \to \infty. \end{cases}$$

REMARK. Note that for $1 < s \leq 2$, as $u \to 0$ and $u/v \to \infty$,

$$v\phi_s(u/v) = v\left\{(1-s) + s\frac{u}{v} - \left(\frac{u}{v}\right)^s\right\}\frac{1}{s(1-s)}$$

$$= \frac{u}{s(s-1)}\left\{\left(\frac{u}{v}\right)^{s-1} + (s-1)\frac{v}{u} - s\right\}$$

$$\sim \frac{u}{s(s-1)}\left\{\left(\frac{u}{v}\right)^{s-1} - s\right\}(1 + o(1)),$$

where the right-hand side converges to $u \log(u/v)(1 + o(1))$ as $s \searrow 1$.

PROOF OF LEMMA 7.2. (i) Letting $u = tv$, it suffices to show that for $0 < v \leq 1/2$ and $1 \leq t \leq 1/(2v)$,

$$K_s(tv, v) \leq \frac{1}{2}\frac{(t-1)^2}{1-v}v$$

or, equivalently, since $K_s(tv, v) = v\phi_s(t) + (1-v)\phi_s((1-tv)/(1-t))$,

$$\phi_s(t) + \left(\frac{1}{v} - 1\right)\phi_s\left(\frac{1-tv}{1-v}\right) \leq \frac{1}{2}\frac{(t-1)^2}{1-v}.$$

Let

$$f_s(t) \equiv \phi_s(t) + \left(\frac{1}{v} - 1\right)\phi_s\left(\frac{1-tv}{1-v}\right) - \frac{1}{2}\frac{(t-1)^2}{1-v}.$$

Now by direct calculation, $f_s(1) = 0$, and

$$f_s'(t) = \phi_s'(t) + \left(\frac{1}{v} - 1\right)\phi_s'\left(\frac{1-tv}{1-v}\right)\left(\frac{-v}{1-v}\right) - \frac{t-1}{1-v}$$

$$= \phi_s'(t) - \phi_s'\left(\frac{1-tv}{1-v}\right) - \frac{t-1}{1-v},$$

so that $f'_s(1) = 0$. Furthermore,

$$f''_s(t) = \phi''_s(t) + \frac{v}{1-v}\phi''_s\left(\frac{1-tv}{1-v}\right) - \frac{1}{1-v}$$

$$= \frac{1}{1-v}\left\{(1-v)\phi''_s(t) + v\phi''_s\left(\frac{1-tv}{1-v}\right) - 1\right\}$$

$$= \frac{1}{1-v}\left\{(1-v)\left(\frac{1}{t}\right)^{2-s} + v\left(\frac{1-v}{1-tv}\right)^{2-s} - 1\right\}$$

$$\leq \frac{1}{1-v}\left\{\left((1-v)\frac{1}{t} + v\left(\frac{1-v}{1-tv}\right)\right)^{2-s} - 1\right\}$$

$$\text{(since } x^{2-s} \text{ is concave for } 1 \leq s \leq 2)$$

$$= \frac{1}{1-v}\left\{\left(\frac{1-v}{t(1-tv)}\right)^{2-s} - 1\right\} \leq 0$$

using the fact that $v(1-v) \leq vt(1-vt)$ for $0 \leq v \leq vt \leq 1/2$ implies $(1-v)/(t(1-tv)) \leq 1$. Here we have used

$$\phi'_s(x) = \frac{s - sx^{s-1}}{s(1-s)} = \frac{1-x^{s-1}}{1-s}, \qquad \phi''_s(x) = x^{s-2}.$$

Since $1 \leq t \leq 1/(2v)$, it follows that $v \leq vt \leq 1/2 < 1$, $1 - v \geq 1 - vt \geq 1/2 > 0$ and $1 \geq (1-vt)/(1-v) \geq 1/(2(1-v))$. When $s < 1$, we calculate

$$f'''_s(t) = \phi'''_s(t) - \left(\frac{v}{1-v}\right)^2 \phi'''_s\left(\frac{1-tv}{1-v}\right)$$

and note that $f'''_s(1) < 0$, while $f'''_s(t) = 0$ has a unique root, so to show $f''_s(t) \leq 0$, it suffices to show $f''_s(1/(2v)) \leq 0$ for $0 \leq v \leq 1/2$. By a straightforward calculation, we get

$$f''_s(1/(2v)) = (2v)^{2-s} + \frac{v}{1-v}\left(\frac{1-v}{1/2}\right)^{2-s} - \frac{1}{1-v},$$

which is $\leq 0$ for $0 \leq v \leq 1/2$ if $s \geq -1$. This shows that $f''_s(t) \leq 0$ in the range $-1 \leq s < 1$, and completes the proof of (i).

(ii) By expanding $K_s(u, v)$ as a function of $u$ as in (10),

$$K_s(u, v) = K_2(u, v)\{1 + v(1-v)D_s(u^*, v) - 1\}$$

with $|u^* - v| \leq |u - v|$; since $0 < v \leq u$, we necessarily have $0 < v \leq u^* \leq u$. Here

$$v(1-v)D_s(u^*, v) - 1 = (1-v)\left\{\left(\frac{v}{u^*}\right)^{2-s} - 1\right\} + v\left\{\left(\frac{1-v}{1-u^*}\right)^{2-s} - 1\right\}$$

$$\equiv I + II.$$

Now $0 < v \leq u^* \leq u$ implies $1 \leq u^*/v \leq u/v$, so $v/u \leq v/u^* \leq 1$ and $1 - v \geq 1 - u^* \geq 1 - u$ implies $(1-v)/(1-u) \geq (1-v)/(1-u^*) \geq 1$. Thus, $I \leq 0$ and $II \geq 0$. It follows that

$$v(1-v)D_s(u^*, v) - 1 \leq v\left\{\left(\frac{1-v}{1-u}\right)^{2-s} - 1\right\}.$$

Similarly,

$$v(1-v)D_s(u^*, v) - 1 \geq (1-v)\left\{\left(\frac{v}{u}\right)^{2-s} - 1\right\}$$

in this range, and the claimed bound in the first part of (ii) follows.

To prove the second part of (ii), note that when $u/v \to \infty$, we can write

$$\frac{K_s(u, v)}{v\phi_s(u/v)} = \frac{1 - (u/v)^s v - ((1-u)/(1-v))^s(1-v)}{v\{1 - s + s(u/v) - (u/v)^s\}}$$

$$= \frac{(u/v)^s v + ((1-u)/(1-v))^s(1-v) - 1}{(u/v)^s v - v(1-s) - su}$$

$$= \frac{1 + [((1-u)/(1-v))^s(1-v) - 1]/[(u/v)^s v]}{1 - [v(1-s) + su]/[(u/v)^s v]}$$

$$\equiv \frac{1 + A(u, v)}{1 - B(u, v)},$$

where, for $1 < s \leq 2$,

$$B(u, v) = \frac{v(1-s) + su}{(u/v)^s v} = \frac{1-s}{(u/v)^s} + s(v/u)^{s-1} = o(1)$$

and

$$A(u, v) = \frac{((1-u)/(1-v))^s(1-v) - 1}{(u/v)^s v}$$

$$= \frac{((1-u)/(1-v))^s[(1-v) - 1] + ((1-u)/(1-v))^s - 1}{(u/v)^s v}$$

$$= -\frac{((1-u)/(1-v))^s}{(u/v)^s} + \frac{((1-u)/(1-v))^s - 1}{(u/v)^s v}$$

$$= o(1) + \frac{1 - su + sv - 1}{(u/v)^s v} + o(1)$$

$$= o(1) + s\frac{1 - (u/v)}{(u/v)^s} = o(1) - s(v/u)^{s-1} = o(1).$$

Thus, the second part of (ii) holds.

The first part of (iv) is proved exactly as in (ii). To prove the second part of (iv), we write

$$K_s(u,v) = \frac{1}{s(1-s)}\{1 - (u/v)^s v - ((1-u)/(1-v))^s(1-v)\}$$

$$= \frac{1}{s(1-s)}\left\{v\left(1 - \left(\frac{u}{v}\right)^s\right) + (1-v)\left(1 - \left(\frac{1-u}{1-v}\right)^s\right)\right\}$$

$$= \frac{1}{s(1-s)}\left\{u\left(\frac{u}{v}\right)^{s-1}\left(-1 + \left(\frac{v}{u}\right)^s\right)\right.$$

$$\left. + \frac{(1-v)(s(u-v) + o(u) + o(v))}{1 - sv + o(v)}\right\}$$

$$= \frac{1}{s(1-s)}\{su(1-(v/u))(1+o(1)) - u(v/u)^{1-s}\{1-(v/u)^s\}\}$$

$$= \frac{1}{1-s}u(1+o(1)).$$

(v) The proof of (v) is similar to the proof of (iv).  $\square$

LEMMA 7.3.  *Suppose that $X_1, \ldots, X_n$ are i.i.d. $F_n$ with $0 < \rho^*(\beta) < r < \beta/3$. Then $r < 1/4$ and for any $0 < r_0 < r$,*

$$(26) \qquad \sup_{n^{-4r} < x < n^{-4r_0}}\left|\frac{\mathbb{F}_n(x)}{x} - 1\right| \xrightarrow{p} 0.$$

PROOF.  Note that $\mathbb{F}_n(\cdot) \overset{d}{=} \mathbb{G}_n(F_n(\cdot))$, where $\mathbb{G}_n$ is the empirical d.f. of $n$ i.i.d. $U(0,1)$ random variables $\xi_1, \ldots, \xi_n$ and

$$F_n(x) = x + \epsilon_n\{(1-x) - \Phi(\Phi^{-1}(1-x) - \mu_n)\} \geq x.$$

Thus, with $\|\cdot\|_a^b \equiv \sup_{a \leq t \leq b}|f(t)|$,

$$\sup_{n^{-4r} < x < n^{-4r_0}}\left|\frac{\mathbb{F}_n(x)}{x} - 1\right|$$

$$(27) \qquad = \left\|\left(\frac{\mathbb{F}_n(x)}{x} - 1\right)^+\right\|_{n^{-4r}}^{n^{-4r_0}} \vee \left\|\left(1 - \frac{\mathbb{F}_n(x)}{x}\right)^+\right\|_{n^{-4r}}^{n^{-4r_0}}$$

$$(28) \qquad \overset{d}{=} \left\|\left(\frac{\mathbb{G}_n(F_n(x))}{x} - 1\right)^+\right\|_{n^{-4r}}^{n^{-4r_0}} \vee \left\|\left(1 - \frac{\mathbb{G}_n(F_n(x))}{x}\right)^+\right\|_{n^{-4r}}^{n^{-4r_0}}.$$

The second term in this last display converges to 0 in probability easily since $F_n(x)/x \geq 1$ implies that it is bounded by

$$\left\|\left(1 - \frac{\mathbb{G}_n(F_n(x))}{F_n(x)}\right)^+\right\|_{n^{-4r}}^1 \leq \left\|\left(1 - \frac{\mathbb{G}_n(t)}{t}\right)^+\right\|_{n^{-4r}}^1 \xrightarrow{p} 0$$

by Theorem 0 of Wellner [42]. On the other hand,

$$\frac{\mathbb{G}_n(F_n(x))}{x} - 1 = \frac{\mathbb{G}_n(F_n(x))}{F_n(x)} \frac{F_n(x)}{x} - 1$$

$$= \left(\frac{\mathbb{G}_n(F_n(x))}{F_n(x)} - 1\right)\frac{F_n(x)}{x} + \left(\frac{F_n(x)}{x} - 1\right),$$

so again by Theorem 0 of Wellner [42], the first term of (27) converges to 0 in probability if

$$\limsup_n \|F_n(x)/x\|_{n^{-4r}}^1 < \infty \quad \text{and} \quad \sup_{n^{-4r} < x < n^{-r_0}} \left(\frac{F_n(x)}{x} - 1\right) \to 0.$$

But this holds by a straightforward analysis using the asymptotics of $\Phi^{-1}$ when $r < \beta/3$. $\quad\square$

Now we have the tools in place to prove our extension of the results of Donoho and Jin [15].

PROOF OF THEOREM 5.1.  First consider $1 < s < 2$. As in Donoho and Jin [15], we first consider the case $r < \beta/3$. Then $r < 1/4$ and we can choose $0 < r_0 < r < 1/4$. From Lemma 7.3 the convergence (26) holds. Thus, by part (ii) of Lemma 7.2, it follows that for $n^{-4r} < x < n^{-4r_0}$, we have

$$nK_s^+(\mathbb{F}_n(x), x) = \frac{1}{2}\left(\frac{(\mathbb{F}_n(x) - x)^+}{\sqrt{x(1-x)}}\right)^2(1 + o_p(1)),$$

and hence,

$$nS_n^+(s) \ge \sup_{n^{-4r} < x < n^{-4r_0}} nK_s^+(\mathbb{F}_n(x), x) \ge \tfrac{1}{2}HC_{n,r,r_0}^{*2}(1 + o_p(1)).$$

Thus, $nS_n^+(s)$ separates $H_0$ and $H_1^{(n)}$ for $s \in (1,2)$ and $r < \beta/3$.

Now suppose that $r > (1 - \sqrt{1-\beta})^2$ (and still $1 < s < 2$). Since $(r + \beta)/(2\sqrt{r}) < 1$, we can pick a constant $q < 1$ such that

$$\frac{(r + \beta)}{2\sqrt{r}} \vee \sqrt{r} < \sqrt{q} < 1.$$

As argued by Donoho and Jin [15], under $H_1^{(n)}$, $n\mathbb{F}_n(n^{-q}) \sim \text{Binomial}(n, L_n n^{-[\beta + (\sqrt{q} - \sqrt{r})^2]})$, where $L_n n^{-[\beta + (\sqrt{q} - \sqrt{r})^2]} \gg n^{-q}$; here $L_n$ is a logarithmic term that does not contribute significantly to the argument. Hence, we have $\mathbb{F}_n(n^{-q})/n^{-q} \gg 1$, and thus, from part (ii) of Lemma 7.2 again,

$$nK_s^+(\mathbb{F}_n(n^{-q}), n^{-q}) = \frac{n\mathbb{F}_n(n^{-q})}{s(s-1)}\left\{\left(\frac{\mathbb{F}_n(n^{-q})}{n^{-q}}\right)^{s-1} - s\right\}(1 + o_p(1)).$$

Hence, we conclude that

$$nS_n^+(s) \geq nK_s^+(\mathbb{F}_n(n^{-q}), n^{-q}) = \frac{n\mathbb{F}_n(n^{-q})}{s(s-1)}\left\{\left(\frac{\mathbb{F}_n(n^{-q})}{n^{-q}}\right)^{s-1} - s\right\}(1 + o_p(1)),$$

so using $\beta + (\sqrt{q} - \sqrt{r})^2 < q < 1$, we conclude that $nS_n^+(s)$ separates $H_0$ and $H_1^{(n)}$ in this range.

Now consider $-1 \leq s < 1$. For this range of $s$ the argument is exactly the same as above, but now using parts (iv) and (v) of Lemma 7.2. (Note that the conclusion of Lemma A.1 of Donoho and Jin [15] can be strengthened considerably as follows: if $Z_n \sim \text{Bin}(n, \pi_n)$ with $\pi_n \to 0$ and $n\pi_n \to \infty$, then $Z_n \overset{p}{\to} \infty$; i.e., for any number $M > 0$, we have $P(Z_n \geq M) \to 1$. This follows easily from Theorem 0 of Wellner [42] since $|Z_n/(n\pi_n) - 1| \overset{p}{\to} 0$ so $Z_n = (Z_n/n\pi_n)n\pi_n \overset{p}{\to} 1 \cdot \infty = \infty$. This also follows easily from the Paley–Zygmund inequality (see, e.g., Kallenberg [29], page 40): $P(Z_n > rE(Z_n)) \geq (1-r)_+^2(EZ_n)^2/[EZ_n^2]$.) □

**Acknowledgments.** We owe thanks to Art Owen for sharing his computer programs, Robert Berk for encouragement, and David Donoho and Jiashun Jin for pointing out several additional references. We also owe thanks to the referees for a number of corrections, suggestions and further work on finite sample computations.

DEPARTMENT OF MATHEMATICS
  AND STATISTICS
GRINNELL COLLEGE
GRINNELL, IOWA 50112-1690
USA
E-MAIL: jagerlea@grinnell.edu

UNIVERSITY OF WASHINGTON
DEPARTMENT OF STATISTICS
BOX 354322
SEATTLE, WASHINGTON 98195-4322
USA
E-MAIL: jaw@stat.washington.edu